\documentclass[]{amsart}
\usepackage{latexsym,amssymb,amsopn,amsbsy,amsthm,amsmath,amscd}

\newcommand{\g}{\mathfrak{g}}
\newcommand{\h}{\mathfrak{h}}
\newcommand{\p}{\mathfrak{p}}

\newcommand{\R}{\mathbb{R}}
\newcommand{\C}{\mathbb{C}}
\newcommand{\Z}{\mathbb{Z}}
\newcommand{\bbS}{\mathbb{S}}

\newcommand{\gC}{\g_{\C}}
\newcommand{\tC}{\mathfrak{t}_{\C}}
\newcommand{\semi}{\ltimes}
\newcommand{\LG}{\tilde{L}G}
\newcommand{\LH}{\tilde{L}H}
\newcommand{\Lg}{\tilde{L}\g}
\newcommand{\Lh}{\tilde{L}\h}
\newcommand{\Waff}{\mathcal{W}}
\newcommand{\blambda}{\boldsymbol{\lambda}}
\newcommand{\bmu}{\boldsymbol{\mu}}
\newcommand{\brho}{\boldsymbol{\rho}}
\newcommand{\balpha}{\boldsymbol{\alpha}}
\newcommand{\Hilbert}{\mathcal{H}}

\newcommand{\SU}{\mathrm{SU}}
\newcommand{\U}{\mathrm{U}}
\newcommand{\SO}{\mathrm{SO}}
\newcommand{\Spin}{\mathrm{Spin}}
\newcommand{\dirac}{\mbox{$\not\negthinspace\partial$}}

\newcommand{\so}{\mathfrak{so}}
\newcommand{\su}{\mathfrak{su}}
\newcommand{\spin}{\mathfrak{spin}}

\DeclareMathOperator{\ad}{ad}

\newcommand{\coad}[1][]{\ad^{\ast}_{#1}}

\DeclareMathOperator{\rank}{rank}
\DeclareMathOperator{\Cl}{Cl}
\DeclareMathOperator{\Diff}{Diff}
\DeclareMathOperator{\tr}{tr}
\DeclareMathOperator{\Sym}{Sym}
\DeclareMathOperator{\ads}
{\widetilde{a\thickspace\thinspace}\negthickspace\negthinspace d}

\DeclareMathOperator{\End}{End}

\DeclareMathOperator{\Ker}{Ker}
\DeclareMathOperator{\Coker}{Coker}
\DeclareMathOperator{\Index}{Index}
\DeclareMathOperator{\sdeg}{sdeg}
\DeclareMathOperator{\ch}{ch}

\theoremstyle{plain}
\newtheorem{theorem}{Theorem}

\newtheorem{lemma}[theorem]{Lemma}
\newtheorem{corollary}[theorem]{Corollary}
\theoremstyle{remark}
\newtheorem*{example}{Example}

\newtheorem*{remark}{Remark}

\begin{document}

\title[Kostant's
Dirac operator for loop groups]{Multiplets of representations and
Kostant's Dirac operator for equal rank loop groups}


\date{October 31, 2000}

\author{Gregory~D.~Landweber}

\address{Microsoft Research\\
         One Microsoft Way\\
         Redmond, WA  98052}

\curraddr{Mathematical Sciences Research Institute\\
         1000 Centennial Drive\\
	 Berkeley, CA  94720-5070}

\email{gregl@msri.org}



\begin{abstract}
Let $\g$ be a semi-simple Lie algebra and let $\h$ be a reductive
subalgebra of maximal rank in $\g$. Given any irreducible
representation of $\g$, consider its tensor product with the spin
representation associated to the orthogonal complement of $\h$ in
$\g$. Gross, Kostant, Ramond, and Sternberg recently proved a
generalization of the Weyl character formula which decomposes the
signed character of this product representation in terms of the
characters of a set of irreducible representations of $\h$, called
a multiplet. Kostant then constructed a formal $\h$-equivariant
Dirac operator on such product representations whose kernel is
precisely the multiplet of $\h$-representations corresponding to
the given representation of $\g$.

We reproduce these results in the Kac-Moody setting for the
extended loop algebras $\Lg$ and $\Lh$. We prove a homogeneous
generalization of the Weyl-Kac character formula, which now yields
a multiplet of irreducible positive energy representations of
$L\h$ associated to any irreducible positive energy representation
of $L\g$. We construct a $L\h$-equivariant operator, analogous to
Kostant's Dirac operator, on the tensor product of a
representation of $L\g$ with the spin representation associated to
the complement of $L\h$ in $L\g$. We then prove that the kernel of
this operator gives the $L\h$-multiplet corresponding to the
original representation of $L\g$.

\end{abstract}

\maketitle

\tableofcontents

\setcounter{section}{-1}

\section{Introduction}

Although this paper is chiefly concerned with representations of
Lie groups and loop groups, the motivation for these results
originally comes from M-Theory. In physics, the Lie group
$\Spin(9)$ arises as the little group for massive particles in 10
dimensional superstring theories and as the little group for
massless particles in 11 dimensional supergravity. Recently,
Pengpan and Ramond noticed that the irreducible representations of
$\Spin(9)$ come in triples, with the Casimir operator taking the
same value on all three representations, and where the dimensions
of two such representations sum to the dimension of the third.
Ramond brought this curious fact to the attention of Sternberg,
who in collaboration with Gross and Kostant then showed that these
triples of representations of $B_{4} = \Spin(9)$ actually
correspond to representations of the exceptional Lie group
$F_{4}$, which contains $B_{4}$ as an equal rank subgroup.

In fact, this is not an isolated phenomenon. In \cite{GKRS},
Gross, Kostant, Ramond, and Sternberg consider the general case
where $\h$ is a reductive Lie algebra which is a maximal rank
subalgebra of some semi-simple Lie algebra $\g$. Letting $G$ and
$H$ denote the compact, simply-connected Lie groups with Lie
algebras $\g$ and $\h$ respectively, associated to any irreducible
representation of $G$ is a set of $\chi(G/H)$ irreducible
representations of $H$, where $\chi(G/H)$ is the Euler number of
the homogeneous space $G/H$. We shall refer to such a set of
$H$-representations as a \emph{multiplet}. As in the case of
$B_{4}\subset F_{4}$, all of the representations in a multiplet
share the same value of the Casimir operator, and the alternating
sum of the dimensions of these representations vanishes. The
relation between a representation of $G$ and the
$H$-representations in the corresponding multiplet is given by the
following homogeneous generalization of the Weyl character
formula, viewed as an identity in the representation ring $R(H)$:
\begin{equation}\label{eq:intro-homogeneous-weyl}
    V_{\lambda} \otimes\bbS^{+} - V_{\lambda} \otimes\bbS^{-}
    = {\sum}_{c\in C}(-1)^{c}\,U_{c(\lambda+\rho_{G})-\rho_{H}},
\end{equation}
where $V_{\lambda}$ and $U_{\mu}$ denote the irreducible
representations of $G$ and $H$ with highest weight $\lambda$ and
$\mu$ respectively, $\bbS = \bbS^{+}\oplus\bbS^{-}$ is the spin
representation associated to the complement of $\h$ in $\g$, the
subset $C\subset W_{G}$ of the Weyl group of $G$ has one
representative from each coset of $W_{H}$, and $(-1)^{c}$ is the
sign of the element $c$.

In representation theory, the Casimir operator of a Lie algebra is
analogous to the Laplacian. Using the spin representation, we can
also consider operators analogous to the Dirac operator.
Furthermore, we can choose a particular Dirac operator such that
its square is the Casimir operator shifted by a constant, giving a
representation theory version of the Weitzenb\"ock formula. Such a
Dirac operator was introduced in a more formal setting by Alekseev
and Meinrenken in \cite{AM}, and the geometric version of this
Dirac operator is examined in \cite{Sl0}. Since the Casimir
operator takes the same value on all of the representations in a
multiplet, it follows that this Dirac operator likewise takes a
constant value, up to sign, on each multiplet.

In the homogeneous case, for any linear operator
\begin{equation*}
    \dirac: V_{\lambda}\otimes\bbS^{+} \rightarrow
    V_{\lambda}\otimes\bbS^{-},
\end{equation*}
since both the domain and range are finite dimensional, the index
of $\dirac$ must be given by (\ref{eq:intro-homogeneous-weyl}).
This prompted Kostant to seach for a Dirac operator whose kernel
and cokernel are precisely those representations on the right hand
side of (\ref{eq:intro-homogeneous-weyl}). In \cite{K}, Kostant
constructs a Dirac operator $\dirac_{\g/\h}$ on
$V_{\lambda}\otimes\bbS$ with a cubic term associated to the
fundamental 3-form on $\g$. The kernel of Kostant's Dirac operator
is
\begin{equation}\label{eq:intro-kernel}
    \Ker\dirac_{\g/\h} = {\bigoplus}_{c\in C}
    U_{c(\lambda+\rho_{G})-\rho_{H}},
\end{equation}
and the signs $(-1)^{c}$ on the right side of
(\ref{eq:intro-homogeneous-weyl}) can be recovered by decomposing
the operator $\dirac_{\g/\h}$ according to the positive and
negative half-spin representations. Taking the kernel of Kostant's
Dirac operator therefore gives an explicit construction of the
multiplet of $H$-representations corresponding to a given
representation of $G$.

This paper takes the results discussed above and reformulates them
in the Kac-Moody setting, replacing the equal rank Lie groups
$H\subset G$ with their corresponding loop groups $LH\subset LG$.
After briefly reviewing the representation theory of loop groups
in \S\ref{section:per}, we introduce the positive energy spin
representation $\mathcal{S}_{L\g}$ associated to a loop group in
\S\ref{section:spin}, using it to reformulate the Weyl-Kac
character formula. In \S\ref{section:weyl-kac}, we prove the
following homogeneous version of the Weyl-Kac character formula:
\begin{equation*}
    \Hilbert_{\blambda} \otimes\mathcal{S}_{L\g/L\h}^{+}
    - \Hilbert_{\blambda} \otimes\mathcal{S}_{L\g/L\h}^{-}
    = {\sum}_{c\in \mathcal{C}}
    (-1)^{c}\,\mathcal{U}_{c(\blambda-\brho_{\g})+\brho_{\h}},
\end{equation*}
where $\Hilbert_{\blambda}$ and $\mathcal{U}_{\bmu}$ denote the
positive energy representations of the central extensions $\LG$
and $\LH$ with lowest weights $\blambda$ and $\bmu$ respectively,
the subset $\mathcal{C}\subset\Waff_{G}$ now lives in the affine
Weyl group of $G$, and $-\brho_{\g}$ and $-\brho_{\h}$ are the
lowest weights of the spin representations $\mathcal{S}_{L\g}$ and
$\mathcal{S}_{L\h}$.

In \S\S\ref{section:clifford}--\ref{section:dirac-gh}, we return
to the case of compact Lie groups, reviewing various results of
\cite{AM} and \cite{K}. There we construct Kostant's Dirac
operator, compute its square, and prove that its kernel has the
form given by (\ref{eq:intro-kernel}). Our approach here differs
slightly from Kostant's, which views the Lie algebra $\g$ as an
orthogonal extension of $\h$. Instead, we first consider the Dirac
operator on $\g$ and a twisted Dirac operator on $\h$ and then
construct Kostant's Dirac operator as their difference, an idea
borrowed from \cite{KS1,KS2}. In addition, we avoid working with a
basis for $\g$ wherever possible, which greatly simplifies the
computations and hopefully elucidates their meanings. These
sections can stand alone as an alternative exposition on Kostant's
Dirac operator, and they provide a outline of the more advanced
material in the subsequent sections.

The remaining sections reprise these results for the loop group
case. In \S\ref{section:loop-clifford} we examine the Clifford
algebra associated to a loop group, which builds on the treatment
of infinite dimensional Clifford algebras given in \cite{KS}. We
then introduce the Dirac and Casimir operators associated to a
loop group in \S\ref{section:dirac-lg}, and we construct the loop
group analogue $\dirac_{L\g/L\h}$ of Kostant's Dirac operator in
\S\ref{section:dirac-lg-lh}. These Dirac and Casimir operators
appear in the physics literature in \cite{KT} and \cite{KS1,KS2}
as the odd and even zero-mode generators for the $N=1$
superconformal algebras associated to current (Lie group) and
coset space (homogeneous space) models. In contrast, our treatment
builds these operators on a mathematical foundation, viewing them
as canonical objects rather than working in terms of a basis.
Finally, we compute the square of the Dirac operator
$\dirac_{L\g/L\h}$, and in \S\ref{section:kernel} we prove that
its kernel is
\begin{equation*}
   \Ker\dirac_{L\g/L\h}
    = {\bigoplus}_{c\in\mathcal{C}}
    \mathcal{U}_{c(\blambda-\brho_{\g})+\brho_{\h}},
\end{equation*}
just as for compact Lie groups. So once again, taking the kernel
of this Dirac operator provides an explicit construction for the
multiplet of representations of $\LH$ corresponding to any given
representation of $\LG$.

\emph{Note.} Anthony Wassermann, who has independently obtained 
results similar to those in this paper, pointed out to me that
with only minor modifications, the arguments presented here provide
a quick proof of the Weyl-Kac character formula.

\section{Loop groups and their representations}

\label{section:per}

\subsection{Loop Groups}

Let $G$ be a compact connected Lie group, and let $LG$ denote the
group of free loops on $G$, i.e., the space of smooth maps from
$S^{1}$ to $G$, where the product of two loops is taken pointwise.
The Lie algebra of the loop group $LG$ is simply the vector space
$L\g$ of loops on the Lie algebra $\g$ of $G$, with brackets again
taken pointwise. The group $\Diff(S^{1})$ of diffeomorphisms of
the circle acts on loop spaces by reparameterizing the loops, and
in particular the subgroup $S^{1}$ of rigid rotations of the
circle acts on $LG$ and $L\g$. This circle action induces a
$\Z$-grading on the complexified Lie algebra $L\gC$, which is the
closure of the direct sum of the Fourier components
$\bigoplus_{k\in\Z}\gC z^{k}$, where $\gC z^{k}$ denotes loops of
the form $z\mapsto X z^{k}$ for $X\in\gC$. We are interested in
those representations of $LG$ which likewise admit a $\Z$-grading
intertwining with the $S^{1}$-action on $LG$, or in other words
representations of the semi-direct product $S^{1}\semi LG$. Such a
representation $\mathcal{E}$ then decomposes into eigenspaces
$\bigoplus_{k\in\Z}\mathcal{E}(k)$ according to the $S^{1}$-weight
$k$, called the \emph{energy}. (This terminology comes from an
analogy with quantum mechanics, where the energies are eigenvalues
of the Hamiltonian operator, which generates time translation.)

\subsection{The central extension}

The representations that we will consider are actually projective
representations of $LG$. To realize them as true representations,
we must introduce a central extension $\LG$ of $LG$ by $S^{1}$.
This is analogous to taking the universal cover of a compact Lie
group, except that here we lift to a circle bundle rather than a
finite cover.

The corresponding central extension of the Lie algebra, which is
called a \emph{Kac-Moody algebra}, is $\Lg = L\g \oplus
\mathbb{R}I$, where $I$ is the infinitesimal generator of the
central $S^{1}$ subgroup. The Lie bracket on the central extension
$\Lg$ is determined by a choice of $\ad$-invariant inner product
on $L\g$. Any $\ad$-invariant inner product on the Lie algebra
$\g$ induces an inner product on the $L\g$ by averaging the
pointwise inner products. For loops $\xi,\eta\in L\g$, this gives
\begin{equation}\label{eq:loop-inner-product}
    \langle \xi, \eta \rangle
    = \frac{1}{2\pi}\int_{0}^{2\pi}
      \bigl\langle \xi(\theta),\eta(\theta)\bigl\rangle\,d\theta,
\end{equation}
which is $\ad$-invariant on $L\g$. To extend this inner product to
$\Lg$, we must actually go one step further and extend it to the
semi-direct sum $\R\mathop{\tilde{\oplus}}\Lg$, where $\R$ is
generated by the infinitesimal rotation $\partial_{\theta}$, and
we define the inner product by
\begin{equation*}
    \langle a\,\partial_{\theta} + \xi + x I,\,
    b\, \partial_{\theta} + \eta + y I \rangle
    = \langle \xi,\eta \rangle - a y - b x
\end{equation*}
for $a,b,x,y\in\R$ and $\xi,\eta\in L\g$. This inner product is
$\ad$-invariant on the extended Lie algebra
$\R\mathop{\tilde{\oplus}}\Lg$ provided that the Lie bracket on
the central extension $\Lg$ is
\begin{equation}\label{eq:central-extension}
    [\xi,\eta]_{\Lg}
    = [\xi,\eta]_{L\g} +
      \langle \xi,\partial_{\theta}\eta\rangle\,I.
\end{equation}
Although this central extension depends on the original choice of
inner product on $\g$, there is a unique $\ad$-invariant inner
product on $\g$ (up to scaling) if $\g$ is \emph{simple}. In this
case, the \emph{universal central extension} corresponds to the
smallest possible scaling for which the Lie algebra $\Lg$
exponentiates to give a central extension $\LG$ of the loop group
$LG$. This smallest inner product on $\g$ is the \emph{basic inner
product}, which is scaled so that the highest root
$\alpha_{\mathrm{max}}$ of $\g$ satisfies
$\|\alpha_{\mathrm{max}}\|^{2} = 2$.

If $G$ is not simple but only semi-simple, then a given projective
representation of $LG$ can still be lifted to a true
representation of some $S^{1}$ extension of $LG$. However, the
universal central extension of $LG$ is no longer a circle bundle,
but rather an extension by the torus $T^{d}$, where $d$ counts the
number of simple components of $G$. At the Lie algebra level, an
$\ad$-invariant inner product on $\g$ can be scaled separately on
each of the simple components, and the central term in the Lie
bracket (\ref{eq:central-extension}) now becomes $d$ separate
terms corresponding to the basic inner products for each of these
components.

\begin{remark}
Let $G$ be simply connected.
Topologically, the invariant inner products on $\g$ correspond to
elements of the Lie algebra cohomology $H^{3}(\g) \cong
H^{3}(G;\R)$ by associating to any inner product its fundamental
3-form $\omega\in\Lambda^{3}(\g^{\ast})$ given by $\omega(X,Y,Z) =
\langle X,[Y,Z] \rangle$ for $X,Y,Z\in\g$. The possible central
extensions of the Lie algebra $L\g$ by a circle thus correspond to
elements of the real cohomology $H^{3}(G;\R)$, and the universal
central extension of $L\g$ is then an extension by the dual space
$K=H_{3}(G;\R)$. On the other hand, the central extensions of the
loop group $LG$ correspond to circle bundles, which are classified
by their Chern classes $c_{1}\in H^{2}(LG;\Z)\cong H^{3}(G;\Z)$ in
the integral lattice of $H^{3}(G;\Z)$. Writing $L=H_{3}(G;\Z)$ for
the dual lattice in $K$, the universal central extension $\LG$ is
an extension of $LG$ by the torus $K/L$. Using the cohomology
spectral sequence for this extension and noting that
$H^{1}(LG;\Z)= H^{2}(G;\Z) = 0$, we obtain the exact sequence
\begin{equation*}
    0\rightarrow H^{1}(\LG;\Z)\rightarrow H^{1}(K/L;\Z)
    \stackrel{d_{2}}{\rightarrow} H^{2}(LG;\Z)
    \rightarrow H^{2}(\LG;\Z)\rightarrow 0.
\end{equation*}
Now, by our construction of the torus $K/L$, we have a canonical
isomorphism $H^{1}(K/L;\Z)\cong H^{3}(G;\Z)$, and we also have a
canonical isomorphism $H^{2}(LG;\Z) \cong H^{3}(G;\Z)$. The map
$d_{2}$ is therefore a homomorphism $d_{2}:H^{3}(G;\Z)\rightarrow
H^{3}(G;\Z)$, and the universality condition becomes the assertion
that $d_{2}$ be the identity map. In particular, if $\LG$ is the
universal central extension, then $d_{2}$ must be an isomorphism,
and it follows that $H^{1}(\LG;\Z)=H^{2}(\LG;\Z)=0$, which in
terms of homotopy implies that $\LG$ is 2-connected. So, whereas
taking the universal cover of a compact Lie group $G$ kills the
obstruction $\pi_{1}(G)$, the loop group $LG$ is already simply
connected, but taking its universal central extension kills the
obstruction $\pi_{2}(LG)$.
\end{remark}

The semi-direct sum $\R\mathop{\tilde{\oplus}}\Lg$ which we
introduced above is the Lie algebra of the semi-direct product
$S^{1}\semi\LG$, and from here on we refer to representations of
$S^{1}\semi\LG$ as representations of $LG$. Given such a
representation, we call the weight of the central $S^{1}$ in $\LG$
the \emph{level} or \emph{central charge}, and since this circle
by definition commutes with the rest of the loop group, it follows
that the level is constant on each irreducible representation.
Unless stated otherwise, from here on we assume that $G$ is simply
connected and simple, we use the basic inner product on $\g$, and
we let $\LG$ denote the universal central extension. However, the
following discussion can be generalized to the semi-simple case by
treating the $d$ simple components separately and viewing the
level as a $d$-vector.

\subsection{Affine roots and the affine Weyl group}

\label{section:affine}

Let $T$ be a maximal torus of $G$. When considering the
representation theory of loop groups, rather than taking the
abelian subgroup $LT$ as the maximal torus of $LG$, we instead use
the maximal torus $S^{1}\times T\times S^{1}$ of $S^{1}\semi \LG$.
Here the first $S^{1}$ factor corresponds to rotation of loops,
while the second comes from the central extension. The Cartan
subalgebra is then $\mathbb{R}\oplus\mathfrak{t}\oplus\mathbb{R}$,
and the weights of $LG$ are of the form $\blambda =
(m,\lambda,h)$, where $m$ is the energy, $\lambda$ is a weight of
$G$, and $h$ is the level. In this notation, the roots of $LG$,
also called the \emph{affine roots} of $G$, consist of the weights
$(m,\alpha,0)$ with $m\in\Z$ and $\alpha$ a root of $G$, as well
as the weights $(m,0,0)$ for nonzero $m$, counted with
multiplicity $\rank G = \dim\mathfrak{t}$. Given a system of
positive roots for $G$, we take the positive roots of $LG$ to be
the roots $(0,\alpha,0)$ for $\alpha>0$, as well as all roots
$(m,\alpha,0)$ with $m>0$, including roots of the form $(m,0,0)$.
If $\{\alpha_{i}\}$ is a set of simple roots for $G$, then the
corresponding simple affine roots for $LG$ are $(0,\alpha_{i},0)$,
as well as the root $(1,-\alpha_{\text{max}},0)$, where
$\alpha_{\text{max}}$ is the highest root of $G$.

The affine Weyl group $\Waff_{G}$ of $G$ is the group generated by
the reflections through the hyperplanes corresponding to the
affine roots of $G$. In terms of loop groups, given any root
$\balpha = (k,\alpha,0)$ of $LG$, there is a corresponding
$\su(2)$ subalgebra of $\Lg$ generated by the loops
$E_{\alpha}z^{k}$ and $E_{-\alpha}z^{-k}$ and the coroot
\begin{equation*}
    H_{k,\alpha}
    =\bigl[E_{\alpha}z^{k},E_{-\alpha}z^{-k}\bigr]_{\Lg}
    =H_{\alpha}+\tfrac{1}{2}ik\|H_{\alpha}\|^{2}I,
\end{equation*}
where $\{E_{\alpha}, E_{-\alpha}, H_{\alpha}\}$ span the $\su(2)$
subalgebra of $\g$ associated to the root $\alpha$. Note that
these elements are normalized so that $\langle
E_{\alpha},E_{-\alpha}\rangle = \frac{1}{2}\|H_{\alpha}\|^{2} =
2\langle\alpha,\alpha\rangle^{-1}$. The reflection of a weight
$\blambda = (m,\lambda,h)$ through the hyperplane orthogonal to
$\balpha$ is then
\begin{equation}\label{eq:reflection}\begin{split}
    s_{k,\alpha}(\blambda) &=
    \blambda-\blambda(H_{k,\alpha})\balpha \\
    &= \bigl(
        m - \lambda(H_{\alpha})k +
        \tfrac{1}{2}h\|H_{\alpha}\|^{2}k^{2},
        \lambda - \lambda(H_{\alpha})\alpha +
        \tfrac{1}{2}h\|H_{\alpha}\|^{2}k\alpha,
        h \bigr).
\end{split}\end{equation}
Furthermore, these $s_{k,\alpha}$ are generated by the reflections
$s_{0,\alpha}$, which act solely on the $\mathfrak{t}^{\ast}$
component and generate the usual Weyl group $W_{G}$, as well as
the transformations
\begin{equation*}
    t_{\alpha} ( \blambda )
    = s_{1,\alpha} s_{0,\alpha} (\blambda)
    = \bigl( m+\lambda(H_{\alpha})+\tfrac{1}{2}h\|H_{\alpha}\|^{2},
        \lambda+ hH_{\alpha},h\bigr),
\end{equation*}
where we use the inner product to identify the coroot
$H_{\alpha}\in\mathfrak{t}$ with the weight
$\frac{1}{2}\|H_{\alpha}\|^{2}\alpha$ in
$\mathfrak{t}^{\ast}$. Restricting to $\mathfrak{t}^{\ast}$,
the $t_{\alpha}$ are simply translations by the coroots,
which generate the coweight lattice $L\subset\mathfrak{t}$. We
therefore have $\Waff_{G}\cong W_{G}\semi L$.

Note that under the action of the affine Weyl group, the level $h$
is fixed, while the energy $m$ is shifted so as to preserve the
inner product
\begin{equation}\label{eq:inner-product}
    (m_{1},\lambda_{1},h_{1}) \cdot (m_{2},\lambda_{2},h_{2})
    = \langle \lambda_{1},\lambda_{2} \rangle
    - m_{1}h_{2} - m_{2}h_{1}
\end{equation}
on $\R\oplus\mathfrak{t}^{\ast}\negthinspace\oplus\R$. Thus, at
any given level $h$, the affine Weyl action is completely
determined by its restriction to $\mathfrak{t}^{\ast}$. In
particular, the element $s_{k,\alpha}$ corresponds to the
reflection through the hyperplane given by the equation
$\langle\lambda,\alpha\rangle=hk$. These hyperplanes divide
$\mathfrak{t}^{\ast}$ into connected components called
\emph{alcoves}, and the affine Weyl group acts simply transitively
on these alcoves. Given a positive root system for $LG$, the
corresponding \emph{fundamental alcove} is the unique alcove
satisfying $\blambda\cdot\balpha \leq 0$ for all $\balpha>0$. This
alcove is bounded by the hyperplanes corresponding to the
negatives of the simple affine roots, or in other words, a weight
$\blambda = (m,\lambda,h)$ lies in the fundamental alcove if and
only if $-\lambda$ is in the positive Weyl chamber for $G$ and
$\langle\lambda,-\alpha_{\text{max}}\rangle \leq h$.

\subsection{Positive energy representations}

A representation $\Hilbert$ of $LG$ is a \emph{positive energy
representation} if $\Hilbert(k) = 0$ for all $k<m$ for some fixed
integer $m$, or in other words, there is a minimum energy when
$\Hilbert$ is decomposed into its constant energy eigenspaces.
In the literature, positive energy
representations are often normalized so that this minimum energy
is 0. However, we will consider positive energy representations
with the full spectrum of minimum energies. When restricted to the
positive energy representations, the representation theory of loop
groups behaves quite analogously to the representation theory of
compact Lie groups. In particular, the positive energy
representations satisfy the following fundamental properties (for
a complete discussion, see \cite{PS}):
\begin{enumerate}
\item[(i)]
A positive energy representation is \emph{completely reducible}
into a direct sum of (possibly infinitely many) irreducible
positive energy representations.
\item[(ii)]
An irreducible positive energy representation $\Hilbert$ is of
\emph{finite type}: each of the constant energy subspaces
$\Hilbert(k)$ is a finite dimensional representation of $G$.
\item[(iii)]
Every irreducible positive energy representation $\Hilbert$ has a
unique \emph{lowest weight} $\blambda = (m,\lambda,h)$, in the
sense that $\blambda - \balpha$ is not a weight of $\Hilbert$ for
any positive root $\balpha$ of $LG$. The lowest weight space is
one dimensional and generates $\Hilbert$.
\item[(iv)]
A weight $\blambda=(m,\lambda,h)$ is \emph{anti-dominant} for $LG$
if it lies in the fundamental Weyl alcove described at the end of
\S\ref{section:affine} above. The lowest weight of a positive
energy representation is anti-dominant, and every anti-dominant
weight is realized as the lowest weight of some positive energy
representation.
\end{enumerate}
As a consequence of (iii), an irreducible positive energy
representation $\Hilbert$ is completely characterized by its
minimum energy $m$, its minimum energy subspace $\Hilbert(m)\cong
V_{-\lambda}$, and its level $h$. Property (iv) implies that for a
positive energy representation, the level $h$ is always
non-negative and is zero only for the trivial representation.
Also, for a fixed minimum energy $m$, there are only finitely many
positive energy representations at each level $h$, but as the
level tends to infinity, the representation theory of $LG$
resembles that of $G$.

If $\Hilbert_{\blambda}$ is the irreducible positive energy
representation with lowest weight $\blambda=(0,\lambda,h)$, then
$\Hilbert_{\blambda}$ also contains all the weights in the orbit
of $\blambda$ under the affine Weyl group $\Waff_{G}$. Recalling
that the affine Weyl group action preserves the inner product
(\ref{eq:inner-product}), it turns out that the orbit of
$\blambda$ consists of all weights $\bmu=(m,\mu,h)$ at level $h$
satisfying $\blambda\cdot\blambda=\bmu\cdot\bmu$, or equivalently
$\|\mu\|^{2} - 2mh = \|\lambda\|^{2}$. This equation sweeps out a
paraboloid, and the weights of $\Hilbert_{\blambda}$ all lie in
its interior. (As the level $h$ tends to infinity, this paraboloid
flattens into a cone.) For an example, see
Figure~\ref{fig:spin-lsu2} at the end of \S\ref{section:weyl-kac},
which gives the weights of the irreducible representation of
$L\SU(2)$ with lowest weight $(0,-1,2)$.

\section{The spin representation}

\label{section:spin}

If $V$ is a finite dimensional vector space with an inner product,
and $V=W\oplus W^{\ast}$ is a polarization of $V$ into a maximal
isotropic subspace $W$ and its dual, then the spin representation
of the Clifford algebra $\Cl(V)$ can be written in the form
\begin{equation}\label{eq:spin-representation}
    \bbS_{V}=\Lambda^{*}(W)\otimes(\det W)^{-\frac{1}{2}},
\end{equation}
where $\det W$ denotes the top exterior power of $W$. The
resulting spin representation $\bbS_{V}$ is independent of
the choice of polarization, which is accounted for by the factor
of $(\det W)^{-1/2}$. On the other hand, if $V$ is infinite
dimensional, then this determinant factor does not make sense,
and so we can no longer use (\ref{eq:spin-representation}) to define
the spin representation. Without this determinant factor to
correct for the choice of polarization, different polarizations
give rise to distinct spin representations.
For a general discussion of infinite dimensional Clifford algebras
and their spin representations, see~\cite{KS}.

For our purposes, consider the Lie algebra $L\g$ with the inner
product (\ref{eq:loop-inner-product}) induced by the basic inner
product on $\g$. If we complexify $L\g$, then the orthogonal
complement of the Cartan subalgebra $\tC$ in $L\gC$ decomposes
into the sum of the positive and negative root spaces, each of
which is isotropic with respect to the inner product on $L\gC$. We
can therefore use this polarization to define a positive energy
spin representation associated to the complement of $\mathfrak{t}$
in $L\g$:
\begin{equation}\label{eq:loop-spin-representation}
    \mathcal{S}_{L\g/\mathfrak{t}}
    := \, \bbS_{\g/\mathfrak{t}} \otimes
      \,\Lambda^{\ast} \Bigl( {\bigoplus}_{k>0} \gC z^{k} \Bigr)
    = \, \bbS_{\g/\mathfrak{t}} \otimes \,
      {\bigotimes}_{k>0} \Lambda^{\ast} \bigl( \gC z^{k} \bigr),
\end{equation}
where we have explicitly factored out the contribution
$\bbS_{\g/\mathfrak{t}}$ coming from the constant loops (or zero
modes). Here, we have used the expression
(\ref{eq:spin-representation}) for the spin representation,
except that we have dropped the portion of the
$(\det W)^{-1/2}$ factor coming from the positive energy modes.
If we were to include that factor, it would contribute an overall
anomalous energy shift of
\begin{equation}\label{eq:energy-shift}
    \Bigl({\prod}_{k>0} z^{k\dim\g}\Bigr)^{-\frac{1}{2}}
    = z^{-\frac{1}{2}\sum_{k>0} k \dim\g}
    = z^{\frac{1}{24}\dim\g},
\end{equation}
where in the last equality we use the Riemann zeta function trick
to write the infinite sum as $\sum_{k>0}k = \zeta(-1) =
-\frac{1}{12}$. Fortunately, by normalizing the spin
representation to have minimum energy 0, we can safely ignore this
factor.

For the moment, we are interested only in the character of the
spin representation. The restriction of the character of
$\mathcal{S}_{L\g/\mathfrak{t}}$ to $S^{1}\times T$ is completely
determined by the description (\ref{eq:loop-spin-representation})
of the spin representation. However, in correcting for the
infinite determinant factor, the spin representation acquires a
nonzero central charge.

\begin{theorem}\label{theorem:spin-lg}
    If $G$ is simple, then the central charge of the spin
    representation $\mathcal{S}_{L\g/\mathfrak{t}}$ is
    the value of the quadratic Casimir operator of $\g$ in the
    adjoint representation:
    \begin{equation*}
        c_{G} = \Delta_{\ad}^{\g}
        = -\frac{1}{2}\,{\sum}_{i}(\ad X_{i})^{2}
        = \langle \rho_{G},\alpha_{\mathrm{max}}\rangle + 1,
    \end{equation*}
    where $\rho_{G}$ is half the sum of the positive roots,
    $\alpha_{\max}$ is the highest root of $G$, and
    $\{X_{i}\}$ is an orthonormal basis for $\g$.
\end{theorem}

\begin{proof}
    To compute the central charge of the spin representation
    $\mathcal{S}_{L\g/\mathfrak{t}}$, we extend it to obtain the
    spin representation associated to the entire Lie algebra $L\g$.
    Since the construction of spin representations is
    multiplicative, we have
    \begin{equation*}
        \mathcal{S}_{L\g} \cong \,
        \bbS_{\mathfrak{t}}\otimes\mathcal{S}_{L\g/\mathfrak{t}}.
    \end{equation*}
    These two spin representations have the same central charge
    since they differ only by the finite dimensional factor
    $\bbS_{\mathfrak{t}}$. However, the extended spin representation
    $\mathcal{S}_{L\g}$ admits an
    action of the full Lie algebra $\Lg$, and in fact,
    $\mathcal{S}_{L\g}$ is the direct sum of
    $\dim\bbS_{\mathfrak{t}}$ copies of an irreducible
    positive energy representation of $L\g$. Examining the
    structure of this representation,
    the first three energy levels of $\mathcal{S}_{L\g}$ are as
    follows:
    \begin{align*}
        \mathcal{S}_{L\g}(0) &= \bbS_{\g},\\
        \mathcal{S}_{L\g}(1) &= \bbS_{\g}\otimes\gC, \\
        \mathcal{S}_{L\g}(2) &= \bbS_{\g}\otimes \gC  \: \oplus \:
        \bbS_{\g}\otimes \Lambda^{2}(\gC).
    \end{align*}
    Letting $\alpha$ denote the highest root of $\g$, and $c$ the
    central charge of $\mathcal{S}_{L\g}$,
    the highest weights of $\mathcal{S}_{L\g}(0)$
    and $\mathcal{S}_{L\g}(1)$ are then $(0,\rho,c)$ and
    $(1,\rho+\alpha,c)$ respectively, while the weight
    $(2,\rho+2\alpha,c)$ is \emph{not}
    present in $\mathcal{S}_{L\g}(2)$. The weights
    $(0,\rho,c)$ and $(1,\rho+\alpha,c)$ thus form a complete
    string of weights for the root $\balpha=(1,\alpha,0)$,
    and so they must be related to each other by the affine Weyl
    element $s_{1,\alpha}$,
    the reflection through the hyperplane orthogonal to $\balpha$.
    By (\ref{eq:reflection}), the difference of these weights is
    $(1,\rho+\alpha,c)-(0,\rho,c)=\balpha =-
(0,\rho,c)(H_{1,\alpha})\balpha$,
    so we obtain
    \begin{equation*}
        -1 = (0,\rho,c) (H_{1,\alpha})
          = \rho(H_{\alpha}) - \tfrac{1}{2}\|H_{\alpha}\|^{2} c
          = \langle\rho,\alpha\rangle - c,
    \end{equation*}
    where $\frac{1}{2}\|H_{\alpha}\|^{2}=1$ and
    $\rho(H_{\alpha})=\langle\rho,\alpha\rangle$
    in the basic inner product since $\alpha$ is the highest root.
    The central charge of the spin representation is thus
    $c = \langle\rho,\alpha\rangle + 1$.

    The quadratic Casimir operator of a Lie algebra does not
    depend on the choice of orthonormal basis, and it commutes
    with the action of the Lie algebra. It therefore acts by a
    constant times the identity on each irreducible representation.
    On the irreducible representation of highest weight $\alpha$,
    the value of the Casimir operator is
    $\frac{1}{2}\|\alpha\|^{2}+\langle\alpha,\rho\rangle$.
    In particular, if $G$ is simple, then the adjoint
    representation is irreducible, and taking $\alpha$ to be
    the highest root of $G$, which satisfies $\|\alpha\|^{2}=2$
    in the basic inner product, we again obtain the value
    $\langle\rho,\alpha\rangle + 1$ as desired.
\end{proof}

We can now compute the character of
$\mathcal{S}_{L\g/\mathfrak{t}}$ directly from the decomposition
(\ref{eq:loop-spin-representation}) and
Theorem~\ref{theorem:spin-lg}. Written in terms of the affine
roots $\balpha = (k,\alpha,0)$, the character is
\begin{equation*}
    \chi(\mathcal{S}_{L\g/\mathfrak{t}})
    = u^{c_{G}}\prod_{\alpha>0}
      \bigl(e^{\frac{i\alpha}{2}}+e^{-\frac{i\alpha}{2}}\bigr)
      \prod_{k>0,\,\alpha}
      \bigl(1 + e^{i\alpha}z^{k}\bigr)
    = e^{-i\brho_{G}}\prod_{\balpha>0}(1 + e^{i\balpha}),
\end{equation*}
where $u$ is a parameter on the central $S^{1}$ extension in
$\LG$, and $\brho_{G} = (0,\rho_{G},-c_{G})$. Here, $-\brho_{G}$
is the lowest weight of $\mathcal{S}_{L\g/\mathfrak{t}}$, which
corresponds to the square root of the determinant in
(\ref{eq:spin-representation}). This weight is the loop group
version of $\rho_{G}$, half the sum of the positive roots of $G$,
which is also characterized by the identity $\rho_{G}(H_{\alpha})
= 1$ for each of the simple roots $\alpha$ of $G$. In the loop
group case, the identity $\brho_{G}(H_{\balpha}) = 1$ must hold as
$\balpha$ ranges over the simple \emph{affine} roots, including
the additional root $(1,-\alpha_{\text{max}},0)$. However, in our
proof of Theorem~\ref{theorem:spin-lg}, the condition
$\brho_{G}(H_{1,-\alpha_{\text{max}}}) = 1$ is the same equation
(up to sign) that we used to compute the central charge $c_{G}$.

The spin representation decomposes as
$\mathcal{S}_{L\g/\mathfrak{t}}=\mathcal{S}_{L\g/\mathfrak{t}}^{+}
\oplus \mathcal{S}_{L\g/\mathfrak{t}}^{-}$ into the sum of two
half-spin representations. In particular, since the complement of
$\mathfrak{t}$ in $\g$ is even dimensional, the zero mode factor
$\bbS_{\g/\mathfrak{t}}$ of $\mathcal{S}_{L\g/\mathfrak{t}}$
splits into half-spin representations, and the exterior algebra in
(\ref{eq:loop-spin-representation}) splits into its even and odd
degree components. The difference of the characters of these
half-spin representations is
\begin{equation}\label{eq:loop-super-character}
    \chi\bigl(\mathcal{S}_{L\g/\mathfrak{t}}^{+}\bigr) -
    \chi\bigl(\mathcal{S}_{L\g/\mathfrak{t}}^{-}\bigr)
    = e^{-i\brho_{G}}{\prod}_{\balpha>0}\bigl(1 - e^{i\balpha}\bigr),
\end{equation}
which can be viewed either as a supertrace on
$\mathcal{S}_{L\g/\mathfrak{t}}$ or as the character of the
virtual representation $\mathcal{S}_{L\g/\mathfrak{t}}^{+}
-\mathcal{S}_{L\g/\mathfrak{t}}^{-}$. Using the notation of spin
representations, the Weyl-Kac character formula becomes

\begin{theorem}[Weyl-Kac Character Formula]
    If $G$ is simply connected and simple, then the character of
    the irreducible positive energy representation
    $\Hilbert_{\blambda}$ of $\LG$ with lowest weight $\blambda$ is
    given by the quotient
    \begin{equation}\label{eq:weyl-kac}
        \chi(\Hilbert_{\blambda}) = \frac
        {\sum_{w\in\Waff_{G}}(-1)^{w}e^{iw(\blambda-\brho_{G})}}
        {\chi\bigl(\mathcal{S}_{L\g/\mathfrak{t}}^{+}\bigr) -
\chi\bigl(\mathcal{S}_{L\g/\mathfrak{t}}^{-}\bigr)},
    \end{equation}
    where $\Waff_{G}$ is the affine Weyl group of $G$ and
    $\brho_{G} = (0,\rho_{G},-c_{G})$.
\end{theorem}

Note that as an immediate consequence of the Weyl-Kac character
formula, if we consider the trivial representation with $\blambda
= 0$, we obtain the identity
\begin{equation*}
    \chi\bigl(\mathcal{S}_{L\g/\mathfrak{t}}^{+}\bigr) -
    \chi\bigl(\mathcal{S}_{L\g/\mathfrak{t}}^{-}\bigr)
    = {\sum}_{w\in\Waff_{G}}(-1)^{w}e^{-iw(\brho_{G})},
\end{equation*}
which gives an alternative expression for the signed character
(\ref{eq:loop-super-character}) of the spin representation
appearing in the denominator of (\ref{eq:weyl-kac}).

\begin{remark}
If $G$ is semi-simple, then we recall that the universal central
extension of $LG$ is an extension not by a circle but rather by
the torus $T^{d}$, where $d$ counts the number of simple
components. In this case the central charge of the spin
representation is the $d$-vector $\mathbf{c}_{G} =
(c_{G_{1}},\ldots,c_{G_{d}})$, where $G_{1},\ldots,G_{d}$ are the
simple components of $G$. If we work with the universal central
extension of $LG$ and define $\brho_{G} =
(0,\rho_{G},-\mathbf{c}_{G})$, then the Weyl-Kac character formula
still holds as written. In fact, using the appropriate universal
central extension, the Weyl-Kac character formula continues to
hold for an arbitrary compact Lie group $G$.
\end{remark}

\section{The homogeneous Weyl-Kac formula}

\label{section:weyl-kac}

Let $\g$ be a compact, semi-simple Lie algebra, and let $\h$ be a
reductive subalgebra of maximal rank in $\g$. In \cite{GKRS},
Gross, Kostant, Ramond, and Sternberg prove a homogeneous
generalization of the Weyl character formula, associating to each
$\g$-representation a set of $\h$-representations with similar
properties, called a multiplet.

\begin{theorem}[Homogeneous Weyl Formula]\label{th:homogeneous-weyl}
    Let $V_{\lambda}$ and $U_{\mu}$ denote the irreducible
    representations of $\g$ and $\h$ with highest weights
    $\lambda$ and $\mu$ respectively. The following identity
    holds in the representation ring $R(\h)$:
    \begin{equation}\label{eq:homogeneous-weyl}
        V_{\lambda}\otimes\bbS^{+}_{\g/\h} -
        V_{\lambda}\otimes\bbS^{-}_{\g/\h}
        = {\sum}_{c\in C}(-1)^{c}\,U_{c(\lambda+\rho_{\g})-\rho_{\h}},
    \end{equation}
    where the sum is taken over the subset $C$
    of elements $c\in W_{\g}$ of the Weyl group of $\g$ for which
    $c(\lambda+\rho_{\g})-\rho_{\h}$ are dominant weights of $\h$.
\end{theorem}

Note that if $\h=\mathfrak{t}$ is a Cartan subalgebra of $\g$,
then $C$ is the full Weyl group $W_{\g}$, and
(\ref{eq:homogeneous-weyl}) becomes the Weyl character formula.
Also note that by stating this result in terms of the Lie algebras
$\h\subset\g$ rather than their corresponding Lie groups $H\subset
G$, we bypass the issue of whether the spin representation
$\bbS_{\g/\h}$ exponentiates to give a true representation of $H$.
Geometrically, this is equivalent to the condition that $G/H$ be a
spin manifold.

Theorem~\ref{th:homogeneous-weyl} has an immediate analogue for
loop groups. The only complication is that simply working at the
level of Lie algebras is no longer sufficient to avoid the
geometric obstruction, which in this case is the condition that
$G/H$ admit a \emph{string structure} (see \cite{M}). Rather, we
must work with the universal central extensions. Given $\g$ and
$\h$ as described above, let $\Lg$ be the universal central
extension of $L\g$, and let $\Lh$ be the restriction of $\Lg$ to
$L\h$. Note that $\Lh$ is not in general the universal central
extension of $L\h$, which we denote by $\hat{L}\h$. Rather, $\Lh$
is a quotient of $\hat{L}\h$. Since $\h$ has the same rank as
$\g$, if $\mathfrak{t}$ is a Cartan subalgebra of $\h$, then it is
likewise a Cartan subalgebra of $\g$. The Cartan subalgebras of
$\R\mathop{\tilde{\oplus}}\hat{L}\h$ and
$\R\mathop{\tilde{\oplus}}\Lg$ are then
$\R\oplus\mathfrak{t}\oplus\R^{d_{\h}}$ and
$\R\oplus\mathfrak{t}\oplus\R^{d_{\g}}$ respectively, where
$d_{\g}$ is the number of simple components of $\g$ and
$d_{\h}\geq d_{\g}$. In other words, we have the commutative
diagram
\begin{equation*}
    \begin{CD}
        \R\oplus\mathfrak{t}\oplus\R^{d_{\h}} @>\text{quotient}>>
        \R\oplus\mathfrak{t}\oplus\R^{d_{\g}} @=
        \R\oplus\mathfrak{t}\oplus\R^{d_{\g}} \\
        @VVV  @VVV  @VVV \\
        \R\mathop{\tilde{\oplus}}\hat{L}\h @>\text{quotient}>>
        \R\mathop{\tilde{\oplus}}\Lh       @>\text{inclusion}>>
        \R\mathop{\tilde{\oplus}}\Lg
    \end{CD}
\end{equation*}
where the vertical arrows are inclusions of Cartan subalgebras.

The weights of $L\h$ and $L\g$ live in the dual spaces to their
Cartan subalgebras, and dual to the quotient map we have an
inclusion
\begin{equation*}
    \R\oplus\mathfrak{t}^{\ast}\oplus\R^{d_{\g}}
    \longrightarrow
    \R\oplus\mathfrak{t}^{\ast}\oplus\R^{d_{\h}}.
\end{equation*}
We may therefore view the weight lattice of $L\g$ as a subset of
the weight lattice of $L\h$. On the other hand, if we ignore the
central extension (i.e., restrict to weights of level 0), then the
weight lattices are identical, and the roots of $L\h$ are a subset
of the roots of $L\g$. Consequently, the affine Weyl group
$\Waff_{\h}$ of $\h$, which is generated by the reflections
through the hyperplanes orthogonal to the roots of $L\h$, is a
subgroup of the affine Weyl group $\Waff_{\g}$ of $\g$.

\begin{theorem}[Homogeneous Weyl-Kac Formula]
    \label{th:homogeneous-weyl-kac}
    Let $\Hilbert_{\blambda}$ and $\mathcal{U}_{\bmu}$ denote the
    irreducible positive energy representations of $\Lg$ and
    $\hat{L}\h$ with lowest weights $\blambda$ and $\bmu$
    respectively. We then have the following identity for virtual
    representations of $\hat{L}\h$:
    \begin{equation}\label{eq:homogeneous-weyl-kac}
        \Hilbert_{\blambda} \otimes \mathcal{S}_{L\g/L\h}^{+} -
        \Hilbert_{\blambda} \otimes \mathcal{S}_{L\g/L\h}^{-}
        = {\sum}_{c\in\mathcal{C}}\,(-1)^{c}\,
      \mathcal{U}_{c(\blambda-\brho_{\g})+\brho_{\h}},
    \end{equation}
    where the sum is taken over the subset $\mathcal{C}$ of
    elements $c\in\Waff_{\g}$ of the affine Weyl group of
    $\g$ for which
    $c(\blambda-\brho_{\g})+\brho_{\h}$ are anti-dominant weights
    of $\hat{L}\h$.
\end{theorem}

\begin{proof}
    We first note that the construction of the spin representation
    is multiplicative, provided that the underlying vector spaces
    are even dimensional. In our case, the positive and negative
    energy subspaces pair off, while for the zero
    modes, the maximal rank condition implies that the complement
    of $\h$ in $\g$ and the complement of $\mathfrak{t}$ in $\h$
    are even dimensional, so we have
    \begin{equation}\label{eq:loop-spin-split}
        \mathcal{S}_{L\g/\mathfrak{t}}^{+} -
        \mathcal{S}_{L\g/\mathfrak{t}}^{-} =
        \bigl( \mathcal{S}_{L\g/L\h}^{+} -
          \mathcal{S}_{L\g/L\h}^{-} \bigr)
        \otimes
        \bigl( \mathcal{S}_{L\h/\mathfrak{t}}^{+} -
          \mathcal{S}_{L\h/\mathfrak{t}}^{-} \bigr).
    \end{equation}
    Applying the Weyl-Kac character formula (\ref{eq:weyl-kac}) to
    the left side of (\ref{eq:homogeneous-weyl-kac}), and factoring
    the Weyl-Kac denominator using (\ref{eq:loop-spin-split}),
    we obtain
    \begin{equation}\label{eq:temp1}
        \chi \bigl( \Hilbert_{\blambda} \otimes
                    \mathcal{S}_{L\g/L\h}^{+} \bigr) -
        \chi \bigl( \Hilbert_{\blambda} \otimes
                    \mathcal{S}_{L\g/L\h}^{-} \bigr)
        = \frac{\sum_{w\in\Waff_{\g}}
	   (-1)^{w}e^{iw(\blambda-\brho_{\g})}}
           {\chi \bigl(\mathcal{S}_{L\h/\mathfrak{t}}^{+} \bigr) -
            \chi \bigl(\mathcal{S}_{L\h/\mathfrak{t}}^{-} \bigr)}.
    \end{equation}
    Recall that the affine Weyl group acts simply transitively on
    the Weyl alcoves. Due to the $\brho_{\g}$ shift, the weight
    $\blambda-\brho_{\g}$ lies in the interior of the fundamental
    Weyl alcove for $\g$, and thus for any $w\in\Waff_{\g}$, the
    weight $w(\blambda-\brho_{\g})$ likewise lies in the interior
    of some Weyl alcove. Furthermore, the Weyl alcoves for $\g$ are
    completely contained inside the Weyl alcoves for $\h$, and
    so there exists a unique element $w'\in\Waff_{\h}$ such that
    $w'w(\blambda-\brho_{\g})$ lies in the interior of
    the fundamental Weyl alcove for $\h$. Shifting back by
    $\brho_{\h}$, we see that the weight
    $w'w(\blambda-\brho_{\g})+\brho_{\h}$ is anti-dominant for
    $\hat{L}\h$. Putting $c = w'w$, we can write $w = (w')^{-1}c$,
    and more generally we have $\Waff_{\g} = \Waff_{\h}\,\mathcal{C}$.
    Using this decomposition to rewrite the numerator on the right
    side of (\ref{eq:temp1}), we have
    \begin{equation*}\begin{split}
        \frac{\sum_{w\in\Waff_{\g}}(-1)^{w}e^{iw(\blambda-\brho_{\g})}}
           {\chi \bigl(\mathcal{S}_{L\h/\mathfrak{t}}^{+} \bigr) -
            \chi \bigl(\mathcal{S}_{L\h/\mathfrak{t}}^{-} \bigr)}
        &= \sum_{c\in\mathcal{C}}(-1)^{c}\,
          \frac{\sum_{w\in\Waff_{\h}}(-1)^{w}e^{iwc(\blambda-\brho_{\g})}}
           {\chi \bigl(\mathcal{S}_{L\h/\mathfrak{t}}^{+} \bigr) -
            \chi \bigl(\mathcal{S}_{L\h/\mathfrak{t}}^{-} \bigr)}\\
        &= \sum_{c\in\mathcal{C}}(-1)^{c}\,
           \chi \bigl( \mathcal{U}_{c\,(\blambda-\brho_{\g})+\brho_{\h}}
\bigr),
    \end{split}\end{equation*}
    where the second line follows by applying the Weyl-Kac character
    formula (\ref{eq:weyl-kac}) for $\hat{L}\h$. This proves the
    character form of the identity (\ref{eq:homogeneous-weyl-kac}).
\end{proof}

The subset $\mathcal{C}\subset\Waff_{\g}$ appearing in
Theorem~\ref{th:homogeneous-weyl-kac} does not depend on the
weight $\blambda$. Rather, it consists of all elements of the
affine Weyl group of $\g$ that map the fundamental Weyl alcove for
$\g$ into the fundamental Weyl alcove for $\h$. Since the affine
Weyl group acts simply transitively on the Weyl alcoves, it
follows that the cardinality of $\mathcal{C}$ is the ratio of the
volumes of the fundamental alcoves for $\g$ and $\h$.
Equivalently, the elements of $\mathcal{C}$ are representatives of
the cosets of $\Waff_{\h}$ in $\Waff_{\g}$, so the cardinality of
$\mathcal{C}$ is the index of $\Waff_{\h}$ in $\Waff_{\g}$. In
particular, the sum appearing in (\ref{eq:homogeneous-weyl-kac})
is finite if and only if $\h$ is semi-simple. In such cases,
$|\mathcal{C}|$ is the index of $W_{\h}$ in $W_{\g}$, which is
also the Euler number of the corresponding homogeneous space
$G/H$. Examples of pairs $\h\subset\g$ with both $\h$ and $\g$
semi-simple include $D_{n}\subset B_{n}$ with $|\mathcal{C}|=2$,
as well as the case $B_{4}\subset F_{4}$ with $|\mathcal{C}|=3$
that prompted \cite{GKRS}. On the other hand, for pairs
$\h\subset\g$ corresponding to \emph{complex} homogeneous spaces
$G/H$, the group $H$ must contain a $\U(1)$ component, and so
(\ref{eq:homogeneous-weyl-kac}) is an infinite sum. We note that
in the physics literature (see \cite{KS1,KS2}), the $N=1$
superconformal coset models on $G/H$ possess an additional $N=2$
symmetry precisely when $\mathcal{C}$ is infinite.

At the other extreme, if $\h=\mathfrak{t}$ is a Cartan subalgebra of
$\g$, then $\brho_{\h}$ vanishes, $\mathcal{C}$ is the full affine
Weyl group $\Waff_{\g}$, and the homogeneous Weyl-Kac formula becomes
\begin{equation}\label{eq:lt-weyl-kac}
   \Hilbert_{\blambda} \otimes \mathcal{S}_{L\g/L\mathfrak{t}}^{+} -
   \Hilbert_{\blambda} \otimes \mathcal{S}_{L\g/L\mathfrak{t}}^{-}
   = {\sum}_{w\in\Waff_{\g}}\,(-1)^{w}\,
   \mathcal{U}_{w(\blambda-\brho_{\g})}.
\end{equation}
This identity is equivalent to the Weyl-Kac character formula
(\ref{eq:weyl-kac}), but it is expressed slightly differently.
Since $\mathfrak{t}$ is abelian, the irreducible positive energy
representation $\mathcal{U}_{\bmu}$ of $\tilde{L}\mathfrak{t}$
takes the particularly simple form
\begin{equation*}
    \mathcal{U}_{\bmu}
    = \Sym^{\ast}\Bigl({\bigoplus}_{k>0}\mathfrak{t}_{\C}z^{k}\Bigr)
    = {\bigotimes}_{k>0}\Sym^{\ast}\left(\mathfrak{t}_{\C}z^{k}\right),
\end{equation*}
where $\Sym^{\ast}$ is the symmetric algebra, and the character of
this representation is
\begin{equation}\label{eq:chi-lt}
    \chi( \mathcal{U}_{\bmu} )
    = e^{i\bmu}\,{\prod}_{k>0}(1 + z^{k} + \cdots
    )^{\dim\mathfrak{t}}
    = e^{i\bmu}\,{\prod}_{k>0}(1 - z^{k})^{-\dim\mathfrak{t}}.
\end{equation}
On the other hand, the signed character of the spin representation
on $L\mathfrak{t}/\mathfrak{t}$ is
\begin{equation}\label{eq:chi-spin-lt}
    \chi\bigl(\mathcal{S}_{L\mathfrak{t}/\mathfrak{t}}^{+}\bigr) -
    \chi\bigl(\mathcal{S}_{L\mathfrak{t}/\mathfrak{t}}^{-}\bigr)
    = {\prod}_{k>0} (1 - z^{k})^{\dim\mathfrak{t}},
\end{equation}
since the product in (\ref{eq:loop-super-character}) is taken over
the positive roots $(k,0,0)$, each counted with multiplicity
$\dim\mathfrak{t}$. In particular, the products in the characters
(\ref{eq:chi-lt}) and (\ref{eq:chi-spin-lt}) cancel each other,
yielding the Weyl-Kac character formula for $LT$:
\begin{equation*}
    \chi\bigl(
    \mathcal{U}_{\bmu} \otimes
    \mathcal{S}_{L\mathfrak{t}/\mathfrak{t}}^{+} \bigr)
    - \chi\bigl(
    \mathcal{U}_{\bmu} \otimes
    \mathcal{S}_{L\mathfrak{t}/\mathfrak{t}}^{-} \bigr)
    = e^{i\bmu}.
\end{equation*}
So, multiplying the formula (\ref{eq:lt-weyl-kac}) by the
character (\ref{eq:chi-spin-lt}), we recover the usual form of the
Weyl-Kac character formula (\ref{eq:weyl-kac}) for $LG$.

\begin{example}
Take $\g=\su(2)$ and let $\h=\mathfrak{u}(1)$ be the Cartan
subalgebra of diagonal elements. In this particular case, we can
use the homogeneous Weyl-Kac formula to explicitly compute the
character of the entire spin representation
$\mathcal{S}_{L\g/L\h}$, not just the difference of the two
half-spin representations. Here we have $\brho_{\g} = ( 0,
\rho_{\g}, -c_{\g}) = ( 0, 1, -2 )$, and so the lowest weight of
the spin representation is $-\brho_{\g}=(0,-1,2)$. The half-spin
representation $\mathcal{S}^{+}$ (resp. $\mathcal{S}^{-}$) is
obtained by acting on a lowest weight vector by an even (resp.
odd) number of Clifford multiplications by the positive generators
$E_{+}$ and $E_{\pm}z^{n}$ for $n>0$ of $L\g/L\h$. Since each of
these generators shifts the $\su(2)$ weight by $\pm 2$, the
$\su(2)$ weights of all elements in $\mathcal{S}^{+}$ must be of
the form $4n-1$, while the weights for $\mathcal{S}^{-}$ are all
of the form $4n+1$. Hence the weights of $\mathcal{S}^{+}$ and
$\mathcal{S}^{-}$ are distinct, and thus there is no cancellation
when we take their difference.

Applying (\ref{eq:lt-weyl-kac}) for the case of the trivial
representation with $\blambda = 0$, we obtain
\begin{align*}
    \mathcal{S}_{L\g/L\h}^{+}
    &= {\sum}_{w\in\Waff_{\g}^{+}}\,\mathcal{U}_{w(0,-1,2)}
     = {\sum}_{n\in\Z}\,\mathcal{U}_{(2n^{2}-n, 4n - 1, 2)}, \\
    \mathcal{S}_{L\g/L\h}^{-}
    &= {\sum}_{w\in\Waff_{\g}^{-}}\,\mathcal{U}_{w(0,-1,2)}
     = {\sum}_{n\in\Z}\,\mathcal{U}_{(2n^{2}+n, 4n + 1, 2)},
\end{align*}
where we have explicitly written out the action of $\Waff_{\su(2)}
\cong \Z_{2}\semi\Z$:
\begin{equation*}
    w_{n}^{\pm}\,(m,\lambda,h)
    = ( m \pm \lambda n + h n^{2}, \pm \lambda + 2 h n, h ).
\end{equation*}
Using (\ref{eq:chi-lt}) for $\chi(\mathcal{U}_{\bmu})$, the
characters of the half-spin representations are
\begin{align*}
  \chi\bigl(\mathcal{S}_{L\g/L\h}^{+}\bigr)(z,w,u)
  &= u^{2} \, {\sum}_{n\in\mathbb{Z}} w^{4n-1} z^{2n^{2}-n}
           \, {\prod}_{k>0}(1-z^{k})^{-1}, \\
  \chi\bigl(\mathcal{S}_{L\g/L\h}^{-}\bigr)(z,w,u)
  &= u^{2} \, {\sum}_{n\in\mathbb{Z}} w^{4n+1} z^{2n^{2}+n}
           \, {\prod}_{k>0}(1-z^{k})^{-1},
\end{align*}
where the powers of $z$, $w$, and $u$ correspond to the energy,
$\su(2)$ weight, and level respectively. Combining these half-spin
representations, the total spin representation has character
\begin{equation*}
  \chi\bigl(\mathcal{S}_{L\g/L\h}\bigr)(z,w,u)
  = u^{2}\,{\sum}_{n\in\mathbb{Z}}w^{2n-1}z^{\frac{1}{2}n(n-1)}
          \,{\prod}_{k>0}(1-z^{k})^{-1}.
\end{equation*}
The orbit of the lowest weight $-\brho_{\g}=(0,-1,2)$ under the
affine Weyl group consists of all weights $(m,\lambda,2)$ with
$\lambda$ odd and $m = \frac{1}{8}(\lambda^{2}-1)$. This equation
sweeps out a parabola, and the remaining weights live inside this
parabola, satisfying $m
> \frac{1}{8}(\lambda^{2}-1)$. The weights of
$\mathcal{S}_{L\g/L\h}$ are shown in Figure~\ref{fig:spin-lsu2},
with the orbit of $-\brho_{\g}$ drawn as open circles. The
multiplicity of any such weight can be derived from
(\ref{eq:chi-lt}) and is given by the number of partitions of
$m-\frac{1}{8}(\lambda^{2}-1)$ into positive integers.

\begin{figure}
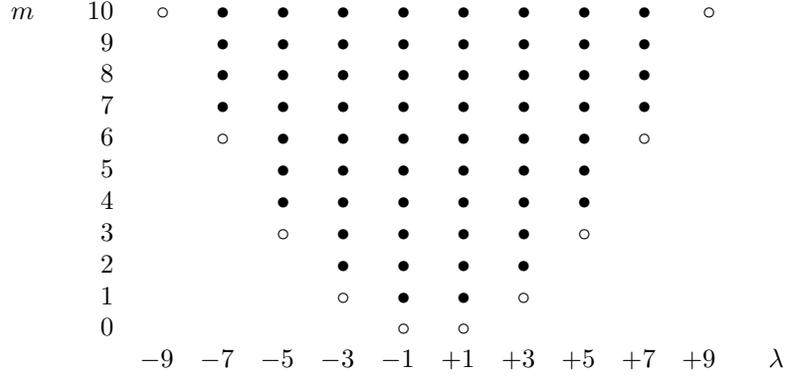
\begin{equation*}\begin{array}{rrrrrrrrrrrr}
m\qquad10&\circ\,&\bullet\,&\bullet\,&\bullet\,&\bullet
\,&\bullet\,&\bullet\,&\bullet\,&\bullet\,&\circ\\
9&&\bullet\,&\bullet\,&\bullet\,&\bullet\,&\bullet
\,&\bullet\,&\bullet\,&\bullet\,&\\
8&&\bullet\,&\bullet\,&\bullet\,&\bullet\,&\bullet
\,&\bullet\,&\bullet\,&\bullet\,&\\
7&&\bullet\,&\bullet\,&\bullet\,&\bullet\,&\bullet
\,&\bullet\,&\bullet\,&\bullet\,&\\
6&&\circ\,&\bullet\,&\bullet\,&\bullet\,&\bullet
\,&\bullet\,&\bullet\,&\circ\,&\\
5&&&\bullet\,&\bullet\,&\bullet\,&\bullet\,&\bullet
\,&\bullet\,&&\\
4&&&\bullet\,&\bullet\,&\bullet\,&\bullet\,&\bullet
\,&\bullet\,&&\\
3&&&\circ\,&\bullet\,&\bullet\,&\bullet\,&\bullet
\,&\circ\,&&\\
2&&&&\bullet\,&\bullet\,&\bullet\,&\bullet\,&&&\\
1&&&&\circ\,&\bullet\,&\bullet\,&\circ\,&&&\\
0&&&&&\circ\,&\circ\,&&&&\\
&-9&-7&-5&-3&-1&+1&+3&+5&+7&+9&\quad\lambda
\end{array}
\end{equation*}
\caption{The weights of the spin representation on
$L\SU(2)/L\U(1)$.}\label{fig:spin-lsu2}\end{figure}

\end{example}
\section{The Clifford algebra $\Cl(\g)$}

\label{section:clifford}

Let $\g$ be a finite dimensional Lie algebra with an
$\ad$-invariant inner product. Recall that the Clifford algebra
$\Cl(\g)$ is generated by the elements of $\g$ subject to the
anti-commutator relation $\{X,Y\} = X\cdot Y + Y\cdot X =
2\,\langle X,Y\rangle$ for all $X,Y\in\g$. There is a natural
Clifford action on the exterior algebra
$\Lambda^{\ast}(\g^{\ast})$, which is given on the generators
$X\in\g$ by $c(X) = \iota_{X} + \varepsilon_{X^{\ast}}$, where
$\iota_{X}$ is interior contraction by $X\in\g$ and
$\varepsilon_{X^{\ast}}$ is exterior multiplication by the dual
element $X^{\ast}\in\g^{\ast}$ satisfying $X^{\ast}(Y)=\langle
X,Y\rangle$. Using the distinguished element $1$ of the exterior
algebra, the map $x \mapsto c(x)1$ gives an isomorphism
$\Cl(\g)\rightarrow\Lambda^{\ast}(\g^{\ast})$ of left
$\Cl(\g)$-modules, called the Chevalley identification. We may
therefore view the Clifford algebra as the exterior algebra
$\Lambda^{\ast}(\g^{\ast})$ with the alternative multiplication
\begin{equation}\label{eq:clifford-multiplication}
    X^{\ast}\!\cdot\eta
    = X^{\ast}\!\wedge\eta + \iota_{X}\eta
\end{equation}
for $X\in\g$ and $\eta\in\Lambda^{\ast}(\g^{\ast})$.

Consider the graded Lie superalgebra $\hat{\g} =
\g_{-1}\oplus\,\g_{0}\oplus\,\R_{1}$, where the subscript denotes
the integer grading. The exterior algebra $\Lambda^{*}(\g^{*})$ is
a representation of this Lie superalgebra $\hat{\g}$, with
$\g_{-1}$ acting by interior contraction, $\g_{0}$ acting by the
coadjoint action, and the generator $d\in\R_{1}$ acting as the
exterior derivative. On the generators $\xi\in\g^{\ast}$, these
operators are given by
\begin{alignat*}{2}
    \iota_{X}\xi &= \xi(X)
    & \qquad\qquad \iota_{X}
    &:\Lambda^{k}(\g^{\ast})\rightarrow\Lambda^{k-1}(\g^{\ast}) \\
    (\coad[X]\xi)(Y) &= -\xi(\ad_{X}Y)
    & \qquad\qquad \ad_{X}^{\ast}
    &:\Lambda^{k}(\g^{\ast})\rightarrow\Lambda^{k+0}(\g^{\ast}) \\
    (d\xi)(X,Y) &= -\tfrac{1}{2}\,\xi([X,Y])
    & \qquad\qquad d
    &:\Lambda^{k}(\g^{\ast})\rightarrow\Lambda^{k+1}(\g^{\ast})
\end{alignat*}
for $X,Y\in\g$. These operators then extend as super-derivations
to the full exterior algebra, and they satisfy the identities
$[\ad_{X}^{\ast},\iota_{Y}] = \iota_{[X,Y]}$ and
$\{d,\iota_{X}\}=\ad_{X}^{\ast}$. If we perturb this action by
taking $d'=d-\iota_{\Omega^{\ast}}$, where $\Omega$ is a closed
$\g$-invariant form of odd degree, then the commutation relations
on $\hat{\g}$ are unchanged.

\begin{theorem}\label{th:adjoint}
    Using the Chevalley identification, the action of
    $\hat{\g}=\g_{-1}\oplus\g_{0}\oplus\R_{1}$ on
    $\Lambda^{\ast}(\g^{\ast})$
    can be expressed in terms of the adjoint action of the Clifford
    algebra as
    \begin{alignat}{2}
    \label{eq:1}
    \ad X^{\ast}
    &= 2\,\iota_{X}
    & \qquad\qquad X^{\ast}   & \in \Lambda^{1}(\g^{\ast}) \\
    \label{eq:2}
    \ad dX^{\ast}
    &= 2\coad[X]
    & \qquad\qquad dX^{\ast} & \in \Lambda^{2}(\g^{\ast}) \\
    \label{eq:3}
    \ad \Omega
    &= 2\,d - 2\,\iota_{\Omega^{\ast}}
    & \qquad\qquad \Omega   & \in \Lambda^{3}(\g^{\ast})
    \end{alignat}
    where $\Omega$ is the fundamental 3-form
    given by $\Omega(X,Y,Z)=-\frac{1}{6}\langle X,[Y,Z]\rangle$.
\end{theorem}

\begin{proof}
    First, we show that the operators $\iota_{X}$ and
    $\ad_{X}^{\ast}$ are super-derivations with respect to
    the Clifford multiplication (\ref{eq:clifford-multiplication}).
    For $X,Y\in\g$ and
    $\eta\in\Lambda^{\ast}(\g^{\ast})$, we have
    \begin{align*}
        \iota_{X}(Y^{\ast}\!\cdot\eta)
        &= \iota_{X}(Y^{\ast}\!\wedge\eta)
        + \iota_{X}\iota_{Y}\eta \\
        &= (\iota_{X}Y^{\ast})\wedge\eta
        - Y^{\ast}\!\wedge\iota_{X}\eta
        - \iota_{Y}\iota_{X}\eta \\
        &= (\iota_{X}Y^{\ast})\cdot\eta
        - Y^{\ast}\!\cdot\iota_{X}\eta,\\
        \coad[X](Y^{\ast}\!\cdot\eta)
        &= \coad[X](Y^{\ast}\!\wedge\eta)
        + \coad[X]\iota_{Y}\,\eta \\
        &= (\coad[X]Y^{\ast})\wedge\eta
        + Y^{\ast}\!\wedge\coad[X]\eta
        + \iota_{Y}\coad[X]\eta
        + \iota_{[X,Y]}\eta \\
        &= (\coad[X]Y^{\ast})\cdot\eta
        + Y^{\ast}\!\cdot \coad[X]\eta.
    \end{align*}
    Now, to prove the identities (\ref{eq:1}) and (\ref{eq:2}), we
    need only verify them for the generators
    $\g^{\ast}=\Lambda^{1}(\g^{\ast})$, but it follows from the
    definition of the Clifford algebra that
    \begin{equation*}
        \{X^{\ast},Y^{\ast}\} = 2\,\langle X,Y\rangle =
        2\,\iota_{X}Y^{\ast},
    \end{equation*}
    and by applying (\ref{eq:1}) and the identity
    $\{d,\iota_{X}\} = \coad[X]$, we obtain
    \begin{equation*}
    [dX^{\ast},Y^{\ast}]
    = -2\,\iota_{Y}dX^{\ast} = -2 \coad[Y]X^{\ast}
    = 2 \coad[X]Y^{\ast}.
    \end{equation*}
    To prove (\ref{eq:3}), we first verify that it holds when
    acting on a generator $X^{\ast}\in\g^{\ast}$:
    \begin{equation*}
    \{\Omega,X^{\ast}\}(Y,Z)
    = ( 2\,\iota_{X}\,\Omega )(Y,Z)
    = -\langle X,[Y,Z]\rangle
    = ( 2\,dX^{\ast})(Y,Z).
    \end{equation*}
    Finally we show that
    $d'=d-\iota_{\Omega^{\ast}}$
    is a super-derivation for Clifford multiplication:
    \begin{equation*}\begin{split}
        d'(X^{*}\!\cdot\eta)
        &= d(X^{*}\!\wedge\eta)
        - \iota_{\Omega^{*}} (X^{\ast}\!\wedge\eta)
        + d\,\iota_{X}\,\eta
        - \iota_{\Omega^{*}} \,\iota_{X}\,\eta \\
        &= (dX^{*})\wedge\eta
        - X^{\ast}\!\wedge d\eta
        - \iota_{(dX^{\ast})^{*}}\,\eta
        + X^{\ast}\!\wedge\iota_{\Omega^{*}}\eta \\
        &\qquad - \iota_{X}d\eta
        + \coad[X]\eta
        + \iota_{X}
            \iota_{\Omega^{*}}\,\eta \\
        &= (d'X^{*})\cdot\eta
        - X^{*}\!\cdot d'\eta,
    \end{split}\end{equation*}
    where we use the expansion
    $
        (d'X^{*})\cdot\eta
        = (dX^{*})\cdot\eta
        = (dX^{*})\wedge\eta
        + \coad[X]\eta
        - \iota_{(dX^{*})^{*}}\,\eta.
    $
\end{proof}

Although the Clifford algebra $\Cl(\g)$ does not admit an integer
grading, it does have the distinguished subspaces $\g$ and
$\spin(\g)$, which correspond via the Chevalley identification to
the first two degrees of the exterior algebra:
\begin{equation*}
    \Lambda^{1}(\g^{\ast})
    \longleftrightarrow \g\subset\Cl(\g),
    \qquad
    \Lambda^{2}(\g^{\ast})
    \longleftrightarrow \spin(\g)\subset\Cl(\g).
\end{equation*}
Since $\Spin(\g)$ is the double cover of $\SO(\g)$, there is a Lie
algebra isomorphism $\spin(\g)\cong\so(\g)$, and given any element
$a\in\so(\g)$, the corresponding element of
$\tilde{a}\in\spin(\g)$ is uniquely determined by the identity
$[\tilde{a},X^{\ast}] = (aX)^{\ast}$ for all $X\in\g$. In
particular, the adjoint action $\ad:\g\rightarrow\so(\g)$ lifts to
a Lie algebra homomorphism $\ads:\g\rightarrow\spin(\g)$
satisfying
\begin{equation*}
    [\ads X,Y^{\ast}] = (\ad_{X}Y)^{\ast} = \coad[X]Y^{\ast}.
\end{equation*}
However, from the identity (\ref{eq:2}), we see that the spin lift
of the adjoint action must be $\ads X = \frac{1}{2}\,dX^{\ast}$.

Let $\{X_{i}\}$ be a basis of $\g$, and let $\{X_{i}^{\ast}\}$
denote the corresponding dual basis of $\g$ satisfying $\langle
X_{i}^{\ast},X_{j}\rangle = \delta_{ij}$. In terms of this basis,
the map $\ads:X\mapsto \frac{1}{2}\,dX^{\ast}$ is
\begin{equation}\label{eq:ads}
    \ads X = -\frac{1}{4}\,{\sum}_{i} X_{i}^{\ast}\cdot [X,X_{i}],
\end{equation}
while the element $\gamma=\frac{1}{4}\,\Omega$ corresponding to
the fundamental 3-form is given by
\begin{equation}\label{eq:gamma}
    \gamma = -\frac{1}{24}\,{\sum}_{i,j}
    X_{i}^{\ast}\cdot X_{j}^{\ast}\cdot [X_{i},X_{j}]
    = \frac{1}{6}\,{\sum}_{i} X_{i}^{\ast} \cdot \ads X_{i}.
\end{equation}
Rewriting Theorem~\ref{th:adjoint} in terms of this new notation,
we obtain the following:

\begin{corollary}\label{cor:superalgebra}
    The elements $1$, $X$, $\ads X$, $\gamma$ for $X\in\g$ span a Lie
    superalgebra $\mbox{$\R_{+}\oplus\g_{-}\oplus\g_{+}\oplus\R_{-}$}$
    in the Clifford algebra $\Cl(\g)$
    with the commutation relations
    \begin{alignat*}{3}
        [\ads X,Y]        &= [X,Y],      & \qquad
        [\ads X, \ads Y]  &= \ads [X,Y], & \qquad
        [\ads X,\gamma]   &= 0, \\
        \{X,Y\}           &= 2\,\langle X,Y\rangle, & \qquad
        \{\gamma,X\}      &= \ads X, & \qquad
        \{\gamma,\gamma\} &= - \tfrac{1}{24}\tr_{\g}\Delta_{\ad}^{\g},
    \end{alignat*}
    where $\Delta_{\ad}^{\g} =
    -\frac{1}{2}\sum_{i}\ad_{X_{i}^{\ast}}\,\ad_{X_{i}}$ is the
    quadratic Casimir operator.
\end{corollary}

\begin{proof}
    All of the commutation relations follow immediately from
    Theorem~\ref{th:adjoint} and the above discussion with the
    exception of that for $\{\gamma,\gamma\}$.  For example, we derive
    \begin{align*}
        [\ads X,\ads Y]
        &= \tfrac{1}{2}\coad[X] dY^{\ast}
        = \tfrac{1}{2}\,d\coad[X] Y^{\ast}
        = \ads[X,Y], \\
        [\ads X,\gamma]
        &= -\bigl[\tfrac{1}{4}\Omega,\tfrac{1}{2}dX^{\ast}\bigr]
        = - \tfrac{1}{4} ddX^{\ast} = 0.
    \end{align*}
    To compute $\{\gamma,\gamma\}$, we note that
    the fundamental 3-form is closed, so we have
    \begin{equation*}\begin{split}
        \{\gamma,\gamma\}
        = \bigl\{\tfrac{1}{4}\,\Omega,\tfrac{1}{4}\,\Omega\bigr\}
        = \tfrac{1}{8}\,d\Omega
        - \tfrac{1}{8}\,\iota_{\Omega^{\ast}}\,\Omega
        = - \tfrac{1}{8}\,\langle\Omega,\Omega\rangle.
    \end{split}\end{equation*}
    Written in terms of an orthonormal basis $\{X_{i}\}$ for $\g$,
    the fundamental 3-form is
    \begin{equation*}
        \Omega =
        -{\sum}_{i<j<k}
        \bigl\langle X_{i},[X_{j},X_{k}]\bigr\rangle \,
        X_{i}^{\ast}\wedge X_{j}^{\ast}\wedge X_{k}^{\ast},
    \end{equation*}
    and so its norm is given by
    \begin{equation*}\begin{split}
        \langle\Omega,\Omega\rangle
        &= \frac{1}{6}\,{\sum}_{i,j,k}
            \bigl\langle X_{i},[X_{j},X_{k}]\bigr\rangle^{2}
        = \frac{1}{6}\,{\sum}_{j,k}
            \bigl\langle [X_{j},X_{k}],[X_{j},X_{k}]\bigr\rangle \\
        &= -\frac{1}{6}\,{\sum}_{j,k}
            \bigl\langle X_{k},[X_{j},[X_{j},X_{k}]]\bigr\rangle
        = \frac{1}{3}\tr_{\g}\Delta_{\ad}^{\g},
    \end{split}\end{equation*}
    which yields the desired anti-commutator
    $\{\gamma,\gamma\}= - \frac{1}{24}\tr_{\g}\Delta_{\ad}^{\g}$.
\end{proof}

Note that the map $\ads:\g\rightarrow\Cl(\g)$ is not necessarily
injective; rather its kernel is the center of $\g$. So,
Corollary~\ref{cor:superalgebra} actually gives us an inclusion of
the superalgebra
\begin{equation}\label{eq:superalgebra}
 \tilde{\g} := \R\oplus\g\oplus[\g,\g]\oplus\langle\g,[\g,\g]\rangle
 \subset\Lambda^{\ast}(\g^{\ast})\cong\Cl(\g)
\end{equation}
into the Clifford algebra of $\g$. Also note that this Lie
superalgebra $\tilde{\g}$, with the commutation relations given by
Corollary~\ref{cor:superalgebra}, is the quantized form of the
graded Lie superalgebra $\hat{\g}=\g_{-1}\oplus\g_{0}\oplus\R_{1}$
discussed above.

\section{The Dirac operator on $\g$}

\label{section:dirac-g}

Let $\g$ be a Lie algebra with an $\ad$-invariant inner product,
and let $\bbS_{\g}$ be the complex spin representation of the
Clifford algebra $\Cl(\g)$. If $\g$ is even dimensional then we
have $\Cl(\g)\otimes\,\C\cong\End(\bbS_{\g})$, and in general the
spin representation $\Cl(\g)\hookrightarrow\End(\bbS_{\g})$ is
faithful. To simplify our notation, in the following we implicitly
identify $\Cl(\g)$ with its image in $\End(\bbS_{\g})$ under the
spin representation. We recall from the previous section that the
adjoint action $\ad$ of $\g$ on itself lifts to the representation
$\ads:\g\rightarrow\Cl(\g)$ given by (\ref{eq:ads}).

Let $V$ be an arbitrary $\g$-module where the $\g$-action is the
map $r:\g\rightarrow\End(V)$. Alternatively, this representation
$r$ may be viewed as the $\End(V)$-valued 1-form
$\hat{r}\in\End(V)\otimes\Lambda^{\ast}(\g^{\ast})$ given
tautologically by $\hat{r}(X) = r(X)$ for all $X\in\g$.
Identifying $\Lambda^{\ast}(\g^{\ast})$ with $\Cl(\g)$ via the
Chevalley map, the element $\hat{r}\in\End(V)\otimes\Cl(\g)$
becomes an operator on the tensor product $V\otimes\bbS_{\g}$.
Perturbing this operator slightly we define the Dirac operator
$\dirac_{r}$ on $V\otimes\bbS_{\g}$ to be the element
\begin{equation*}
    \dirac_{r} := \hat{r} + 1\otimes \tfrac{1}{2}\,\Omega
    \,\in\,\End(V)\otimes\Cl(\g),
\end{equation*}
where $\Omega\in\Cl(\g)$ is the cubic term given by (\ref{eq:3}).
Written in terms of a basis $\{X_{i}\}$ for $\g$ and the dual
basis $\{X_{i}^{\ast}\}$ satisfying $\langle
X_{i}^{\ast},X_{j}\rangle = \delta_{i,j}$, this Dirac operator is
\begin{equation*}\begin{split}
    \dirac_{r}
    &= \sum_{i} X_{i}^{\ast}\,r(X_{i})
    - \frac{1}{12}\sum_{i,j}
        X_{i}^{\ast} \cdot X_{j}^{\ast} \cdot [X_{i},X_{j}] \\
    &= \sum_{i} X_{i}^{\ast}\,\Bigl(
        r(X_{i}) + \frac{1}{3} \ads X_{i} \Bigr).
\end{split}\end{equation*}
Note that the second form of this operator resembles a geometric
Dirac operator for the connection $\nabla_{X}=r(X)+\frac{1}{3}\ads
X$. Indeed, if $r$ is the right action of $\g$ on functions, then
this is the \emph{reductive connection} on the spin bundle over
$G$ (see \cite{Sl0}).

Rather than choosing a particular representation $V$, we can
instead take $r$ to be the canonical inclusion
$r:\g\hookrightarrow U(\g)$ of $\g$ into its universal enveloping
algebra $U(\g)$. This gives us a universal Dirac operator
$\dirac$, which is an element of the non-abelian Weil algebra
$U(\g)\otimes\Cl(\g)$, introduced by Alekseev and Meinrenken in
\cite{AM}. Again identifying $\Cl(\g)$ with
$\Lambda^{\ast}(\g^{\ast})$, the element $\dirac$ is characterized
by the identity
\begin{equation}\label{eq:ix-dirac}
    \iota_{X}\dirac = \varrho(X) := r(X)\otimes 1 + 1\otimes\ads X
\end{equation}
for all $X\in\g$, where $\varrho = r\otimes 1 + 1\otimes\ads$ is
the diagonal action of $\g$ on $U(\g)\otimes\Cl(\g)$. Squaring the
Dirac operator, we obtain the Weitzenb\"ock formula
\begin{equation}\label{eq:dirac-g-square}\begin{split}
    \dirac^{2}
    &= \hat{r}\cdot\hat{r}
       + \bigl\{\hat{r},\tfrac{1}{2}\,\Omega\bigr\}
       + \tfrac{1}{2}\,\Omega \cdot \tfrac{1}{2}\,\Omega \\
    &= \langle\hat{r},\hat{r}\rangle
     + \hat{r}\wedge\hat{r}
     + d\hat{r}
     + \tfrac{1}{8}\,\{ \Omega, \Omega \}
     = -2\,\Delta_{r}^{\g} - \tfrac{1}{12}\tr_{\g}\Delta_{\ad}^{\g},
\end{split}\end{equation}
where the ``curvature'' term $d\hat{r} + \hat{r}\wedge\hat{r}$
vanishes since
\begin{equation*}
    (d\hat{r} + \hat{r}\wedge\hat{r})(X,Y)
    = \tfrac{1}{2} \bigl( - r([X,Y]) + [r(X),r(Y)]\,\bigr)
    = 0.
\end{equation*}
Note that the square of the Dirac operator has no Clifford algebra
component and is thus an element $\dirac^{2}\in U(\g)$ of the
universal enveloping algebra. In fact, since the Casimir operator
commutes with the $\g$-action, the element $\dirac^{2}$ lies in
the center of $U(\g)$. Considering the Dirac operator itself,
given any $X\in\g$ we have
\begin{equation}\label{eq:commutes}
    [\varrho(X),\dirac]
    = [\iota_{X}\dirac,\dirac]
    = \iota_{X}(\dirac\cdot\dirac) = 0,
\end{equation}
and thus $\dirac$ is invariant under the diagonal action $\varrho$
of $\g$ on $U(\g)\otimes\Cl(\g)$. We can summarize the above
results by stating that the Lie superalgebra
\begin{equation*}
    \R \, \oplus \, (1\otimes\g) \, \oplus \, \varrho(\g)
    \, \oplus\, \R\dirac \, \oplus \, \R\Delta^{\g}
    \subset U(\g)\otimes\Cl(\g)
\end{equation*}
is a central extension of the Lie superalgebra $\tilde{\g}$ from
Corollary~\ref{cor:superalgebra} by the span of the quadratic
Casimir operator $\Delta^{\g}$. The commutation relations in this
extension are the same as in Corollary~\ref{cor:superalgebra},
with the exception of the square of the Dirac operator which is
given by (\ref{eq:dirac-g-square}). To obtain the corresponding
``classical'' algebra, we let this superalgebra act on the
non-abelian Weil algebra via the adjoint action. Since the
elements $1$ and $\Delta^{\g}$ lie in the center of the universal
enveloping algebra, we are left with the graded Lie superalgebra
$\hat{\g}$, as Alekseev and Meinrenken show in \cite{AM}.

\begin{theorem}
    The non-abelian Weil algebra $U(\g)\otimes\Cl(\g)$ is a
    representation of the graded Lie superalgebra $\hat{\g} =
    \g_{-1}\oplus\g_{0}\oplus\R_{1}$ spanned by the operators
    $\iota_{X}$, $\mathcal{L}_{X}$, $d$ for $X\in\g$ given by
    \begin{equation*}
        \iota_{X} = \ad \bigl( \tfrac{1}{2}X \bigr),
        \qquad
        \mathcal{L}_{X} = \ad \bigl( \varrho(X) \bigr),
        \qquad
        d = \ad \bigl( \dirac \bigr).
    \end{equation*}
\end{theorem}

Now, suppose that $\g$ is reductive. If $V_{\lambda}$ is the
irreducible representation of $\g$ with highest weight $\lambda$,
then the value of the quadratic Casimir operator
$\Delta_{\lambda}^{\g}$ on $V_{\lambda}$ is
\begin{equation*}
    \Delta_{\lambda}^{\g} = \tfrac{1}{2} \bigl(
    \| \lambda + \rho_{\g} \|^{2} - \| \rho_{\g} \|^{2} \bigr).
\end{equation*}
In addition, for reductive Lie algebras we have the identity
$\frac{1}{12}\tr_{\g}\Delta_{\ad}^{\g}=\|\rho_{\g}\|^{2}$, and it
follows that the square of the Dirac operator $\dirac_{\lambda}$
acting on $V_{\lambda}\otimes\bbS_{\g}$ is simply the constant
$\dirac_{\lambda}^{2} = -\| \lambda + \rho_{\g} \| ^{2}$.

\section{The Dirac operator on $\g/\h$}

\label{section:dirac-gh}

Let $\h$ be a Lie subalgebra of $\g$, and let $\p$ be the
orthogonal complement of $\h$ with respect to the $\ad$-invariant
inner product on $\g$. The adjoint action of $\h$ on
$\g=\h\oplus\p$ respects this decomposition, so we obtain
$\h$-representations $\ad_{\h}$ and $\ad_{\p}$ on $\h$ and $\p$
respectively.
The Clifford algebra also decomposes into the product
$\Cl(\g)\cong\Cl(\h)\otimes\Cl(\p)$ of two Clifford algebras, and
with it the spin lift of the adjoint action becomes the sum
$\ads_{\g} = \ads_{\h}\otimes 1 + 1\otimes\ads_{\p}$ of 
separate spin actions $\ads_{\h}:\h\rightarrow\Cl(\h)$ and
$\ads_{\p}:\h\rightarrow\Cl(\p)$. The spin representations
$\bbS_{\h}$ and $\bbS_{\p}$ of $\Cl(\h)$ and $\Cl(\p)$ are
therefore representations of $\h$, and if one or both of $\h$ or
$\p$ is even dimensional, then we have
$\bbS_{\g}\cong\bbS_{\h}\otimes\bbS_{\p}$.

Let $\dirac_{\g}\in U(\g)\otimes\Cl(\g)$ denote the universal
Dirac operator on $\g$ discussed in the previous section. Now
consider the twisted Dirac operator $\dirac_{\h}'$ on $\h$ given
by
\begin{equation*}
    \dirac_{\h}' := \hat{r}'_{\h} + \tfrac{1}{2}\,\Omega_{\h}
    \in\bigl( U(\h)\otimes\Cl(\p) \bigr) \otimes\Cl(\h)
    \cong U(\h)\otimes\Cl(\g)
\end{equation*}
where $r' = r\otimes 1 + 1\otimes\ads_{\p}$ is the diagonal action
of $\h$ on $U(\h)\otimes\Cl(\p)$. In other words, given any
representation $U$ of $\h$, this twisted Dirac operator
$\dirac_{\h}'$ is the usual Dirac operator $\dirac_{\h}$ acting on
the twisted space $(U\otimes\bbS_{\p})\otimes\bbS_{\h}\cong
U\otimes\bbS_{\g}$. As we saw in (\ref{eq:ix-dirac}), this Dirac
operator $\dirac_{\h}'$ is characterized by the identity
\begin{equation*}
    \iota_{Z}\dirac'_{\h} = \varrho'_{\h}(Z) = \varrho_{\g}(Z)
    = \iota_{Z}\dirac_{\g}
\end{equation*}
for all $Z\in\h$, where $\varrho'_{\h}$ is the diagonal action of
$\h$ on $\bigl( U(\h)\otimes\Cl(\p) \bigr) \otimes\Cl(\h)$. Note
that $\varrho'_{\h}$ is just the restriction to $\h$ of the
diagonal action $\varrho_{\g}$ of $\g$ on $U(\g)\otimes\Cl(\g)$.
It then follows from (\ref{eq:commutes}) that the element
$\dirac_{\h}'$ commutes with the diagonal action $\varrho'_{\h}$.

Define $\dirac_{\g/\h} = \dirac_{\g} - \dirac_{\h}'$ to be the
difference of these two operators. This element $\dirac_{\g/\h}$
is then \emph{basic} with respect to $\hat{\h}$, or in other words
it satisfies the identities
\begin{equation*}
    \iota_{Z}\,\dirac_{\g/\h} = 0, \qquad
    \mathcal{L}_{Z}\,\dirac_{\g/\h}
    = \bigl[\varrho'_{h}(Z),\dirac_{\g/\h}\bigr] = 0,
\end{equation*}
for all $Z\in\h$. This Dirac operator can also be written as the
element
\begin{equation*}
    \dirac_{\g/\h} =
    \hat{r}_{\p} + 1\otimes\tfrac{1}{2}\,\Omega_{\p}
    \in \bigl( U(\g) \otimes \Cl(\p) \bigr)^{\h},
\end{equation*}
where $\hat{r}_{\p}\in U(\g)\otimes\Lambda^{1}(\p^{\ast})$
corresponds to the map $r:\p\hookrightarrow U(\g)$, and
$\Omega_{\p}\in\Lambda^{3}(\p^{\ast})$ is the fundamental 3-form
given by $\Omega_{\p}(X,Y,Z)=-\frac{1}{6}\langle X,[Y,Z]\rangle$
for all $X,Y,Z\in\p$. To keep track of the cubic terms, note that
the fundamental 3-form decomposes as
\begin{equation*}
    \Omega_{\g}
    = \Omega_{\h} + \Omega_{\p} + 3\,\Omega_{\h\p\p},
\end{equation*}
into its projections onto $\Lambda^{3}(\h^{\ast})$,
$\Lambda^{3}(\p^{\ast})$, and
$\h^{\ast}\wedge\p^{\ast}\wedge\p^{\ast}$ respectively. This extra
contribution corresponds to the twisted term $\hat{r}'_{\h} -
\hat{r}_{\h} = \tfrac{3}{2}\,\Omega_{\h\p\p}$ appearing in
$\dirac_{\h}'$.

If $\{X_{i}\}$ and $\{X_{i}^{\ast}\}$ are dual bases for $\p$,
then the Dirac operator $\dirac_{\g/\h}$ is
\begin{equation*}\begin{split}
    \dirac_{\g/\h}
    &= \sum_{i} X_{i}^{\ast}\,r(X_{i})
    - \frac{1}{12}\sum_{i,j}
        X_{i}^{\ast} \cdot X_{j}^{\ast} \cdot [X_{i},X_{j}]_{\p} \\
    &= \sum_{i} X_{i}^{\ast}\,\Bigl(
        r(X_{i}) + \frac{1}{3} \ads_{\p} X_{i} \Bigr),
\end{split}\end{equation*}
where $[X,Y]_{\p}$ for $X,Y\in\p$ denotes the projection of
$[X,Y]$ onto $\p$, and
\begin{equation*}
    \ads_{\p} X
    = -\frac{1}{4}\,{\sum}_{i} X_{i}^{\ast}\cdot [X,X_{i}]_{\p}
\end{equation*}
for $X\in\p$. (Note that all of the sums here are taken over a
basis of $\p$, not of $\g$.) The geometric version of this Dirac
operator $\dirac_{\g/\h}$, viewed as an operator on twisted
spinors on the homogeneous space $G/H$, is discussed in
\cite{Sl1,Sl2} and \cite{L}.

To compute the square of $\dirac_{\g/\h}$, we first show that
$\dirac_{\h}'$ and $\dirac_{\g/\h}$ decouple,
\begin{equation*}
    \bigl\{ \dirac_{\h}',\dirac_{\g/\h} \bigr\}
    = \bigl\{ \hat{r}'_{\h}, \dirac_{\g/\h} \bigr\}
    + \bigl\{ \tfrac{1}{2}\,\Omega_{\h}, \dirac_{\g/\h} \bigr\}
    = [r'(\cdot),\dirac_{\g/\h} ]\tilde{\ }+d_{\h}\dirac_{\g/\h}
    = 0.
\end{equation*}
We therefore have
\begin{equation*}
    \dirac_{\g/\h}^{2}
    = (\dirac_{\g})^{2} - (\dirac_{\h}')^{2}
    = - 2 \bigl( \Delta^{\g}_{r} - \Delta^{\h}_{r'} \bigr)
      - \tfrac{1}{12} \bigl( \tr_{\g}\Delta^{\g}_{\ad}
                           - \tr_{\h}\Delta^{\h}_{\ad} \bigr).
\end{equation*}
Suppose that both $\g$ and $\h$ are reductive. If
$r:\g\rightarrow\End(V_{\lambda})$ is the irreducible
representation of $\g$ with highest weight $\lambda$, then
$\dirac_{\g/\h}$ is an $\h$-invariant operator on
$V_{\lambda}\otimes\bbS_{\p}$. Its square then takes the value
\begin{equation}\label{eq:dirac-gh-square}
    \dirac_{\g/\h}^{2}|_{\mu}
    = - \| \lambda + \rho_{\g} \|^{2}
      + \| \mu + \rho_{\h} \|^{2}
\end{equation}
on the $\h$-invariant subspace of $V_{\lambda}\otimes\bbS_{\p}$
transforming like the irreducible representation $U_{\mu}$ of $\h$
with highest weight $\mu$. It follows that the kernel of
$\dirac_{\g/\h}^{2}$, which is in turn the kernel of the Dirac
operator $\dirac_{\g/\h}$ itself, consists of all $\h$-invariant
subspaces of $V_{\lambda}\otimes\bbS_{\p}$ transforming like
$U_{\mu}$, where $\mu$ is a dominant weight of $\h$ satisfying
$\|\mu + \rho_{\h}\|^{2} = \|\lambda + \rho_{\g}\|^{2}$. As we
show in the following theorem, these subspaces are precisely the
multiplet of $\h$-representations corresponding to the
$\g$-representation $V_{\lambda}$, which we discussed in
\S\ref{section:weyl-kac}.

\begin{theorem}\label{theorem:dirac-gh-kernel}
    Let $\g$ be a semi-simple Lie algebra with a maximal rank
    reductive Lie subalgebra $\h$, and let
    $V_{\lambda}$ and $U_{\mu}$ denote the irreducible
    representations of $\g$ and $\h$ with highest weights
    $\lambda$ and $\mu$. The kernel of the Dirac
    operator $\dirac_{\g/\h}$ on
    $V_{\lambda}\otimes\bbS_{\p}$ is
    \begin{equation*}
        \Ker\dirac_{\g/\h}
        = {\bigoplus}_{c\in C} U_{c\bullet\lambda},
    \end{equation*}
    where $c\bullet\lambda = c(\lambda + \rho_{\g})-\rho_{\h}$, and
    $C\subset W_{\g}$ is the subset of Weyl elements which map the
    positive Weyl chamber for $\g$ into the positive Weyl chamber
    for $\h$.
\end{theorem}

\begin{proof}
    Since the Weyl group acts by isometries, the weights
    $c\bullet\lambda$ satisfy the identity
    \begin{equation*}
        \|(c\bullet\lambda)+\rho_{\h}\|^{2}
        = \|c(\lambda+\rho_{\g})\|^{2}
        = \|\lambda + \rho_{\g}\|^{2},
    \end{equation*}
    and it follows from (\ref{eq:dirac-gh-square}) that
    the Dirac operator $\dirac_{\g/\h}$
    vanishes on any $\h$-invariant subspace of
    $V_{\lambda}\otimes\bbS_{\p}$ transforming like
    $U_{c\bullet\lambda}$. To complete the proof, it remains
    to show that each of these representations occurs exactly
    once in the domain of the Dirac operator and that no other
    $\h$-representations appear in its kernel. We establish
    these facts in the following two lemmas (see also \cite{K0}).
\end{proof}

\begin{lemma}\label{lemma:1}
For each $c\in C$, the irreducible representation
$U_{c\bullet\lambda}$ of $\h$ with highest weight $c\bullet\lambda
= c(\lambda+\rho_{\g}) - \rho_{\h}$ occurs exactly once in the
decomposition of $V_{\lambda}\otimes\bbS_{\p}$.
\end{lemma}

\begin{proof}
The highest weight space of an irreducible representation of $\g$
is always one dimensional, so the weight $\lambda$ appears with
multiplicity 1 in $V_{\lambda}$. Now consider the complex spin
representation $\bbS_{\g/\mathfrak{t}}$ associated to the
orthogonal complement of a Cartan subalgebra $\mathfrak{t}$ in
$\g$. Given a positive root system for $\g$, the character of this
spin representation is $\chi(\bbS_{\g/\mathfrak{t}}) =
\prod_{\alpha>0}(e^{i\alpha/2} + e^{-i\alpha/2})$, and so the
highest weight $\rho_{\g} = \frac{1}{2}\sum_{\alpha>0}\alpha$ of
$\bbS_{\g/\mathfrak{t}}$ appears with multiplicity 1. The highest
weight of the tensor product
$V_{\lambda}\otimes\bbS_{\g/\mathfrak{t}}$ is then $\lambda +
\rho_{\g}$, appearing with multiplicity 1, and likewise the
weights $w(\lambda+\rho_{\g})$ for $w\in W_{\g}$ all have
multiplicity 1 in $V_{\lambda}\otimes\bbS_{\g/\mathfrak{t}}$.

Choosing a common Cartan subalgebra
$\mathfrak{t}\subset\mathfrak{h}\subset{\g}$, the spin
representation factors as
$\bbS_{\g/\mathfrak{t}}\cong\bbS_{\p}\otimes\bbS_{\h/\mathfrak{t}}$.
As we noted above, the weight $\rho_{\h}$ appears with
multiplicity 1 in the second factor $\bbS_{\h/\mathfrak{t}}$. It
follows that the weights $w\bullet\lambda$ for $w\in W_{\g}$ can
appear at most once in the tensor product
$V_{\lambda}\otimes\bbS_{\p}$, as each such weight contributes one
weight of the form $(w\bullet\lambda)+\rho_{\h} =
w(\lambda+\rho_{\g})$ to the tensor product
$V_{\lambda}\otimes\bbS_{\p}\otimes\bbS_{\h/\mathfrak{t}}\cong
V_{\lambda}\otimes\bbS_{\g/\mathfrak{t}}$. On the other hand, we
see from the homogeneous Weyl formula (\ref{eq:homogeneous-weyl})
that the irreducible representations $U_{c\bullet\lambda}$ for
$c\in C$ appear at least once in the decomposition of the tensor
product $V_{\lambda}\otimes\bbS_{\g/\h}$. We therefore conclude
that the representations $U_{c\bullet\lambda}$ for $c\in C$ each
occur exactly once in $V_{\lambda}\otimes\bbS_{\p}$.
\end{proof}

\begin{lemma}\label{lemma:2}
If $\mu$ is a weight of $V_{\lambda}\otimes\bbS_{\p}$ satisfying
$\| \mu+\rho_{\h} \|^{2} = \| \lambda+\rho_{\g} \|^{2}$, then
there exists a unique Weyl element $w\in W_{\g}$ such that
$\mu+\rho_{\h} = w(\lambda+\rho_{\g})$.
\end{lemma}

\begin{proof}
If $\mu$ is a weight of $V_{\lambda}\otimes\bbS_{\p}$, then
$\mu+\rho_{\h}$ is a weight of the tensor product
$V_{\lambda}\otimes\bbS_{\p}\otimes\bbS_{\h/\mathfrak{t}}\cong
V_{\lambda}\otimes\bbS_{\g/\mathfrak{t}}$. Since the Weyl group
acts simply transitively on the Weyl chambers, there exists an
element $w\in W_{\g}$ such that $w^{-1}(\mu+\rho_{\h})$ is
dominant, where we recall that a weight $\nu$ is dominant if and
only if $\langle\nu,\alpha\rangle \geq 0$ for all positive roots
$\alpha$. Note that every weight of $\bbS_{\g/\mathfrak{t}}$ can
be obtained from its highest weight $\rho_{\g}$ by subtracting a
sum of positive roots. Likewise, for the tensor product
$V_{\lambda}\otimes\bbS_{\g/\mathfrak{t}}$, the difference
$(\lambda + \rho_{\g}) - w^{-1}(\mu+\rho_{\h})$ is a sum of
positive roots, and it follows that
\begin{equation*}
    \| \lambda +\rho_{\g} \|^{2}
    \geq \| w^{-1}(\mu+\rho_{\h}) \|^{2},
\end{equation*}
with equality holding only when $(\lambda + \rho_{\g}) -
w^{-1}(\mu+\rho_{\h}) = 0$. As for the uniqueness of $w$, if
$\lambda$ is dominant, then the weight $\lambda + \rho_{\g}$ lies
in the interior of the positive Weyl chamber for $\g$, and thus
the weights $w(\lambda+\rho_{\g})$ for $w\in W_{\g}$ are distinct.
\end{proof}

Theorem~\ref{theorem:dirac-gh-kernel} now follows immediately from
the above two lemmas. We can actually be slightly more specific
about the kernel of the Dirac operator, recovering the signs
appearing in the homogeneous Weyl formula
(\ref{eq:homogeneous-weyl}). Recall that the spin representation
decomposes as  $\bbS_{\p} = \bbS_{\p}^{+}\oplus\bbS_{\p}^{-}$ into
two half-spin representations. Since the Dirac operator is an odd
element of the non-abelian Weil algebra $U(\g)\otimes\Cl(\g)$, it
interchanges $\bbS_{\p}^{+}$ and $\bbS_{\p}^{-}$. Restricting the
domain of the Dirac operator to the positive half-spin
representation, we obtain an operator
\begin{equation*}
    \dirac_{\g/\h}^{+} : V_{\lambda}\otimes\bbS_{\p}^{+}
    \rightarrow V_{\lambda}\otimes\bbS_{\p}^{-}.
\end{equation*}
Furthermore, since the Dirac operator is formally self-adjoint,
its adjoint is
\begin{equation*}
    \dirac_{\g/\h}^{-} : V_{\lambda}\otimes\bbS_{\p}^{-}
    \rightarrow V_{\lambda}\otimes\bbS_{\p}^{+},
\end{equation*}
the restriction of the Dirac operator to the negative half-spin
representation. Since these Dirac operator are acting on finite
dimensional vector spaces, the index is the difference of the
domain and range, so we have
\begin{equation}\label{eq:dirac-index}
    \Ker\dirac_{\g/\h}^{+} - \Ker\dirac_{\g/\h}^{-}
    = V_{\lambda}\otimes\bbS_{\p}^{+}
    - V_{\lambda}\otimes\bbS_{\p}^{-},
\end{equation}
which is given by the homogeneous Weyl formula
(\ref{eq:homogeneous-weyl}). Comparing this with the kernel of
$\dirac_{\g/\h} = \dirac_{\g/\h}^{+} \,\oplus\,
\dirac_{\g/\h}^{-}$ given in
Theorem~\ref{theorem:dirac-gh-kernel}, we therefore obtain
\begin{equation*}
    \Ker\dirac_{\g/\h}^{+} = \bigoplus_{(-1)^{c} = +1}
    U_{c\bullet\lambda}, \qquad
    \Ker\dirac_{\g/\h}^{-} = \bigoplus_{(-1)^{c} = -1}
    U_{c\bullet\lambda}.
\end{equation*}
In other words, there is no cancellation on the left hand side of
equation (\ref{eq:dirac-index}), and the signed kernel of this
Dirac operator picks out precisely those representations, with
sign, appearing on the right hand side of the homogeneous Weyl
formula (\ref{eq:homogeneous-weyl}).

\section{The Clifford algebra $\Cl(L\g)$}

\label{section:loop-clifford}

In Section~\ref{section:clifford},
we examined the Clifford algebra associated to a finite
dimensional Lie algebra with an invariant inner product. The
infinite dimensional case is more complicated, and the
general theory of such infinite dimensional Clifford algebras and
their spin representations is developed in the mathematical
literature by Kostant and Sternberg in \cite{KS}. Here, we
consider the Clifford algebra associated to the Lie algebra $L\g$
of smooth maps from $S^{1}$ to a finite dimensional Lie algebra
$\g$, where we restrict to the dense subspace of loops with finite
Fourier expansions. This finiteness condition ensures that the Lie
algebra $L\g$ has a countable basis, and the complexification of
this loop space is then $L\gC=\bigoplus_{k\in\Z}\gC
z^{k}=\gC[z,z^{-1}]$, the Lie algebra of finite Laurent series
with values in $\gC$. Averaging the pointwise inner products over
the loop, we obtain an invariant inner product on $L\g$ given by
(\ref{eq:loop-inner-product}).

The Clifford algebra $\Cl(L\g)$ is spanned by \emph{finite} sums of
products of the form $\xi_{1}\cdots\xi_{n}$ for loops $\xi_{i}\in
L\g$, subject to the relation $\{\xi,\eta\} =
2\langle\xi,\eta\rangle$.  However, the loop space analogues of the
elements $\ads X$ and $\gamma$ introduced in \S\ref{section:clifford}
are in fact infinite sums, so we must instead work with a formal
completion of the Clifford algebra.  Unfortunately, the product of two
such infinite sums does not necessarily converge.  On the other hand,
given a spin representation $\mathcal{S}_{L\g}$ of the Clifford
algebra $\Cl(L\g)$, we can view $\End(\mathcal{S}_{L\g})$ as a
completion of $\Cl(L\g)$ with a well defined product given by the
composition of endomorphisms.  As we discussed in
\S\ref{section:spin}, to define the spin representation we must first
choose a polarization.  With respect to the action of the
infinitesimal generator $\partial_{\theta}$ of rotations, the
complexified loop space $L\gC$ decomposes into its negative, zero, and
positive energy subspaces $L\gC=L\gC^{-}\oplus\gC\oplus L\gC^{+}$,
where $L\gC^{+}$ and $L\gC^{-}$ are isotropic subspaces which are dual
to each other with respect to the inner product.  The spin
representation corresponding to this polarization is
$\mathcal{S}_{L\g}:=\bbS_{\g}\otimes\Lambda^{*}(L\gC^{+})$, and the
Clifford action $c:\Cl(L\gC)\rightarrow\End(\mathcal{S}_{L\g})$ is
given by
\begin{equation*}\label{eq:loop-clifford-action}
    c(\xi) = \begin{cases}
        1 \otimes \varepsilon(\xi) & \text{for $\xi\in L\gC^{+}$}, \\
        1 \otimes \iota(\xi) & \text{for $\xi\in L\gC^{-}$}, \\
        c(\xi) \otimes (-1)^{F} & \text{for $\xi\in\gC$},
    \end{cases}
\end{equation*}
where $\varepsilon$ and $\iota$ are exterior multiplication and
interior contraction respectively, and $F$ is the degree operator on
the exterior algebra.  If $\{\eta_{i}\}$ is a basis for
$\Cl(L\g^{-}_{\C})$, then when applied to a specific element of the
spin representation $\mathcal{S}_{L\g}$, all but finitely many of the
operators $c(\eta_{i}) = \iota(\eta_{i})$ vanish.  Formal infinite
sums $\sum_{i}c(\omega_{i})c(\eta_{i})$ with coefficients
$\omega_{i}\in\Cl(\gC\oplus L\g^{+}_{\C})$ therefore yield well
defined operators on the spin representation, and in fact all elements
of $\End(\mathcal{S}_{L\g})$ can be expressed in this form.

The exterior algebra that we shall consider here is not
$\Lambda^{*}(L\g^{*})$, but rather the algebra $\Lambda^{*}(L\g)^{*}$
of skew-symmetric multilinear forms on $L\g$.  Such forms can be
expressed as formal infinite sums of basic products of the form
$\xi_{1}^{*}\wedge\cdots\wedge\xi_{n}^{*}$ for $\xi_{i}\in
L\g$.  In infinite dimensions, the Chevalley map
$\ch:\Cl(L\g)\rightarrow\Lambda^{*}(L\g)^{*}$ is no longer surjective;
its image consists of all
forms given by finite sums of the basic wedge products.  Although the
Chevalley map fails to converge if we attempt to extend it to the
completion $\End(\mathcal{S}_{L\g})$ of $\Cl(L\g)$, we can perturb it
by terms of lower degree to remove the infinite contributions. 
Separating the Clifford algebra into its positive and negative energy
factors, we define the \emph{normal ordering} map
$n:\Cl(L\gC)\rightarrow\Lambda^{*}(L\gC)^{*}$ by
\begin{equation*}
    n(\omega^{+}\cdot\omega^{-})
    = \ch(\omega^{+}) \wedge \ch(\omega^{-}),
\end{equation*}
where $\omega^{+}\in\Cl(\gC\oplus
L\g^{+}_{\C})$ and $\omega^{-}\in\Cl(L\g^{-}_{\C})$.
The normal ordering map extends to the completion 
$\End(\mathcal{S}_{L\g})$ of the Clifford algebra, and its image
is the subspace
$\Lambda^{*}(L\gC)^{+}\subset\Lambda^{*}(L\gC)^{*}$ given by
\begin{equation*}\label{eq:normal-ordering}
    \Lambda^{*}(L\gC)^{+}
    = \bigl\{ \omega\in\Lambda^{*}(L\gC)^{*} \bigm|
    (\iota_{\eta}\omega)^{+}
    \in \Lambda^{*}(L\gC^{-})
    \text{ for all } \eta\in\Lambda^{*}(L\gC^{-}) \bigr\},
\end{equation*}
where $(\,)^{+}$ denotes the projection of $\Lambda^{*}(L\gC)^{*}$
onto $\Lambda^{*}(L\gC^{+})^{*}$, and we identify
$\Lambda^{*}(L\gC^{-})$ with a subspace of $\Lambda^{*}(L\gC^{+})^{*}$
via the inner product.  In terms of a basis $\{\eta_{i}\}$ for
$\Lambda^{*}(L\g^{-})$, we may write elements of
$\Lambda^{*}(L\gC)^{+}$ as formal infinite sums
$\sum_{i}\omega_{i}^{*}\wedge\eta_{i}^{*}$, with
$\omega_{i}\in\Lambda^{*}(\gC\oplus L\gC^{+})$ living in the zero and
positive energy components.

\begin{remark}
    Decomposing $L\gC=\bigoplus_{k\in\Z}\gC z^{k}$ in terms of its
    energy grading, we define a secondary degree on $L\gC$
    counting only the negative energy contribution
    \begin{equation*}
	\sdeg Xz^{k} =
	\begin{cases}
	    0  & \text{for $k\geq 0$}, \\
	    k & \text{for $k < 0$},
	\end{cases}
    \end{equation*}
    where $X\in\gC$ and $Xz^{k}$ is the loop $z\mapsto Xz^{k}$ for
    $|z|=1$.  Let $L\gC^{*} = \gC^{*}[z,z^{-1}]$ denote the reduced
    dual of $L\gC$.
    Extending $\sdeg$ to the exterior algebra
    $\Lambda^{*}(L\gC^{*})$, we note that 
    $\Lambda^{*}(L\gC)^{+}$ is the completion of
    $\Lambda^{*}(L\gC^{*})$ with respect to $\sdeg$.  In other words,
    $\Lambda^{*}(L\gC)^{+}$ consists of all formal infinite sums
    $\sum_{i}\omega_{i}$ of $\sdeg$-homogeneous elements
    $\omega_{i}\in\Lambda^{*}(L\gC^{*})$ with
    $\sdeg\omega_{i}\rightarrow\infty$.
\end{remark}

We can now use the normal ordering identification
$n:\End(\mathcal{S}_{L\g})\rightarrow\Lambda^{*}(L\gC)^{+}$ to define
product and bracket structures on $\Lambda^{*}(L\gC)^{+}$.  The normal
ordered product $\omega_{1}\cdot_{n}\omega_{2} =
n(n^{-1}\omega_{1}\cdot n^{-1}\omega_{2} )$ on the exterior algebra
differs from the product induced by the Chevalley identification by
terms of lower degree.  However, many of the supercommutators remain
unchanged.  In particular, the normal ordered bracket with the
dual $\xi^{*}\in L\gC^{*}\cong\Lambda^{1}(L\gC)^{+}$ of a loop $\xi\in
L\gC$ is still given by
\begin{equation*}\begin{split}
    [\xi^{*},\omega^{+}\wedge\omega^{-}]_{n}
    &= n\bigl[n^{-1}(\xi^{*}),n^{-1}(\omega^{+}\wedge\omega^{-})\bigr]
    = n[\xi,\ch^{-1}\omega^{+}\cdot\ch^{-1}\omega^{-}] \\
    &= n\bigl( [\xi,\ch^{-1}\omega^{+}]\cdot\ch^{-1}\omega^{-}
       \pm \ch^{-1}\omega^{+} \cdot [\xi,\ch^{-1}\omega^{-}] \bigr) \\
    &= 2\iota_{\xi}\omega^{+} \wedge \omega^{-}
       \pm \omega^{+}\wedge 2\iota_{\xi}\omega^{-}
     = 2\,\iota_{\xi}(\omega^{+}\wedge\omega^{-}),
\end{split}\end{equation*}
for $\omega^{+}\in\Lambda^{*}(\gC\oplus L\g^{+}_{\C})^{*}$ and
$\omega^{-}\in\Lambda^{*}(L\g^{-}_{\C})^{*}$ of homogeneous
degree.  Thus,
\begin{equation}\label{eq:interior1}
    [ \xi^{*}, \omega ]_{n} = 2\,\iota_{\xi}\,\omega
    \text{ for $\xi\in L\gC$ and $\omega\in\Lambda^{*}(L\gC)^{+}$.}
\end{equation}
Reprising the discussion of \S\ref{section:clifford},
for any $\xi\in L\g$, consider the 2-form $d\xi^{*}$ given by
$d\xi^{*}(\eta,\zeta) =
-\frac{1}{2}\langle\xi,[\eta,\zeta]\rangle$ for all $\eta,\zeta\in
L\g$. Although $\partial_{\theta}$ is not an element of $L\g$, we
can nevertheless define an analogous 2-form
$d\partial_{\theta}^{*}$ by $d\partial_{\theta}^{*}(\xi,\eta) :=
\frac{1}{2} \langle\xi,\partial_{\theta}\eta\rangle$. Note that
$d\partial_{\theta}^{*}$ is closed but not exact, so it defines a
cohomology element in $H^{2}(L\g)$. Finally, the fundamental
3-form $\Omega$ is given by $\Omega(\xi,\eta,\zeta) =
-\frac{1}{6}\langle\xi,[\eta,\zeta]\rangle$ for $\xi,\eta,\zeta\in
L\g$. These elements all lie in $\Lambda^{*}(L\gC)^{+}$, and they
satisfy the identities
\begin{equation}\label{eq:interior2}
    \iota_{\xi}\,d\eta^{\ast}  = [\xi,\eta]^{*},
    \qquad
    \iota_{\xi}\,d\partial_{\theta}^{*} =
    -(\partial_{\theta}\xi)^{*},
    \qquad
    \iota_{\xi}\,\Omega = d\xi^{*}.
\end{equation}
Using the normal ordered product and bracket on
$\Lambda^{*}(L\gC)^{+}$ coming from $\End(\mathcal{S}_{L\g})$, we
obtain the loop space version of Corollary~\ref{cor:superalgebra}.

\begin{theorem}\label{th:loop-superalgebra}
If $\g$ is simple, then the elements $1$, $\xi^{*}$ for $\xi\in
L\g$, $\ads\xi = \frac{1}{2}\,d\xi^{*}$ for $\xi\in
\R\mathop{\tilde{\oplus}}L\g$, and $\gamma=\frac{1}{4}\,\Omega$
span a Lie superalgebra in
$\Lambda^{*}(L\g)^{+}\subset\End(\mathcal{S}_{L\g})$ satisfying
\begin{alignat*}{2}
    \{\xi^{*},\eta^{*}\} &= 2\langle\xi,\eta\rangle,
    \vphantom{\tfrac{1}{2}}\\
    [\ads\xi,\eta^{*}]  &= [\xi,\eta]^{*},
    &\qquad
    [\ads\xi,\ads \eta] &= \ads\,[\xi,\eta] + i c_{\g}
    \langle\xi,\partial_{\theta}\eta\rangle,
    \vphantom{\tfrac{1}{2}}\\
    [\ads\partial_{\theta},\xi^{\ast}] &= (\partial_{\theta}\xi)^{*},
    &\qquad
    [\ads\partial_{\theta},\ads\xi] &= \ads(\partial_{\theta}\xi),
    \vphantom{\tfrac{1}{2}}\\
    \{ \gamma, \xi^{\ast} \} &= \ads\xi,
    &\qquad
    [\gamma,\ads\xi] &= \tfrac{1}{2}\,ic_{\g}(\partial_{\theta}\xi)^{*},
    \\
    [\gamma,\ads\partial_{\theta}] &= 0,
    &\qquad
    \{ \gamma, \gamma \} &=
    i c_{\g}\ads\partial_{\theta} - \tfrac{1}{24}\,c_{\g}\dim\g,
\end{alignat*}
where $c_{\g}$ is the value of the Casimir operator of $\g$ in the
adjoint representation.
\end{theorem}

\begin{proof}
    The bracket $\{\xi^{*},\eta^{*}\} = 2\langle\xi,\eta\rangle$
    is simply the definition of the Clifford algebra, while the
    brackets $[\ads\xi,\eta^{*}] = [\xi,\eta]^{*}$ and
    $[\ads\partial_{\theta},\xi^{\ast}] =
    (\partial_{\theta}\xi)^{*}$ and $\{ \gamma,\xi^{\ast} \} =
    \ads\xi$ follow immediately from (\ref{eq:interior1}) and
    (\ref{eq:interior2}). By the
    Jacobi identity, for any $\xi,\eta\in\R\mathop{\tilde{\oplus}}L\g$
    and $\zeta\in L\g$ we have
    \begin{equation*}\begin{split}
        \bigl[[\ads\xi,\ads\eta],\zeta^{*}\bigr]
        &= \bigl[\ads\xi,[\ads\eta,\zeta^{*}]\bigr]
        - \bigl[\ads\eta,[\ads\xi,\zeta^{*}]\bigr] \\
        &=\bigl[\xi,[\eta,\zeta]\bigr]^{*}-\bigl[\eta,[\xi,\zeta]\bigr]^{*}
        = \bigl[[\xi,\eta],\zeta\bigr]^{*}
        = \bigl[\ads\,[\xi,\eta],\zeta^{*}\bigr],
    \end{split}\end{equation*}
    which shows that $\ads$ is a projective representation of
    $\R\mathop{\tilde{\oplus}}L\g$ on $\mathcal{S}_{L\g}$.
    In Theorem~\ref{theorem:spin-lg}, we established that this
    spin representation has central charge $c_{\g}$, which gives
    us the brackets
    $[\ads\xi,\ads \eta] = \ads\,[\xi,\eta] + ic_{\g}
    \langle\xi,\partial_{\theta}\eta\rangle$ and
    $[\ads\partial_{\theta},\ads\xi] = \ads\partial_{\theta}\xi$.

    To compute $\gamma^{2}$, we write it as the sum $\gamma^{2} =
    (\gamma^{2})_{0} + (\gamma^{2})_{2}$ of homogeneous forms of
    degrees 0 and 2. (We shall see that $\gamma^{2}$ has no
    components of degrees 4 or 6.) Since $\iota_{\xi}=\ad \xi^{*}$ for
    $\xi\in L\g$ is a derivation with respect to
    the backet, we have
    \begin{equation*}
        \iota_{\xi}\gamma^{2}
        = [\iota_{\xi}\gamma,\gamma]
        = \tfrac{1}{2}\,[\ads\xi,\gamma].
    \end{equation*}
    Taking one further interior contraction, we obtain
    \begin{equation*}
        \iota_{\xi}\iota_{\eta}\gamma^{2}
        = - \tfrac{1}{4}\bigl(
        [\ads\xi,\ads\eta] - \ads\,[\xi,\eta] \bigr)
        = - \tfrac{1}{4}\,
        i c_{\g}\langle\xi,\partial_{\theta}\eta\rangle,
    \end{equation*}
    which is a constant. It follows that $\gamma^{2}$ has no
    components of degree higher than 2, and that $(\gamma^{2})_{2}$
    is the 2-cocycle determining the central extension of $L\g$ for
    the spin representation $\ads$. In fact, this 2-cocycle is
    a multiple of $\ads\partial_{\theta}$, and we have
    \begin{equation*}
        (\gamma^{2})_{2}(\xi,\eta)
        = -\tfrac{1}{2}\,\iota_{\xi}\iota_{\eta}\gamma^{2}
        = \tfrac{1}{8}\bigl( i c_{\g} \langle
        \xi,\partial_{\theta}\eta\rangle \bigr)
        = \tfrac{1}{2}\,i c_{\g} (\ads\partial_{\theta})(\xi,\eta).
    \end{equation*}
    Going back up one level, we see that
    \begin{equation*}
        [\ads\xi,\gamma] = 2\,\iota_{\xi}\gamma^{2}
        = i c_{\g}\,\iota_{\xi}\ads\partial_{\theta}
        = -\tfrac{1}{2}\,i c_{\g} (\partial_{\theta}\xi)^{*}.
    \end{equation*}
    Finally, the value of the constant $(\gamma^{2})_{0}$ is
    the value of $\gamma^{2}$ acting on the minimum energy
    subspace $\mathcal{S}_{L\g}(0)$ of the spin representation,
    since $\ads\partial_{\theta}$ vanishes there. However, all
    the terms in $\gamma^{2}$ vanish on $\mathcal{S}_{L\g}(0)$
    except the contribution from the constant loops,
    and thus $(\gamma^{2})_{0} =
    \tfrac{1}{48}\tr_{\g}\Delta_{\ad}^{\g}
    = \tfrac{1}{48}\,c_{\g}\dim\g$
    as we proved in Corollary~\ref{cor:superalgebra}.
\end{proof}

Taking a slightly different view of this theorem, the commutation
relations given in Theorem~\ref{th:loop-superalgebra} determine a
Lie superalgebra (with subscripts denoting the grading)
\begin{equation*}
    \R_{\text{even}} \, \oplus \, L\g_{\text{odd}} \, \oplus \,
    (\R\mathop{\tilde{\oplus}}L\g)_{\text{even}}\,\oplus\,\R_{\text{odd}},
\end{equation*}
and the identification of $\Lambda^{*}(L\g)^{+}$ with its image in
$\End(\mathcal{S}_{L\g})$ gives a representation of this Lie
superalgebra on the spin representation $\mathcal{S}_{L\g}$. 
Actually, we can extend this Lie superalgebra further.  The component
$\R_{\text{even}} \oplus L\g_{\text{odd}} \oplus L\g_{\text{even}}$ is
called a \emph{super Kac-Moody algebra}, and using superspace
notation, its complexification is a central extension of the
polynomial algebra $\g\otimes\C[z,z^{-1},\Theta]$, where $\Theta$ is
an odd variable (i.e. $\Theta^{2}=0$).  The \emph{super Virasoro
algebra} $\mathrm{SVir}$ is the universal central extension of the Lie
algebra of derivations of $\C[z,z^{-1},\Theta]$.  (Note that the even
derivations are just the vector fields on the circle.)  The super
Virasoro algebra therefore acts on the super Kac-Moody algebra, and
their semidirect sum is referred to as the \emph{$N=1$ superconformal
current algebra} (see \cite{KT}):
\begin{equation*}
    \mathrm{SVir} \; \tilde{\oplus} \;
    (\R_{\text{even}}\oplus L\g_{\text{odd}}\oplus L\g_{\text{even}}).
\end{equation*}
In our case, the elements $\ads \partial_{\theta}$ and $\gamma$
span the even and odd zero-mode subspaces of the super Virasoro
algebra, with commutator $\{\gamma,\gamma\} = ic_{\g} \bigl(
\ads\partial_{\theta} + \tfrac{1}{24}\,i\dim\g \bigr)$. Here, the
additional $\tfrac{1}{24}\dim\g$ term, which is sometimes
incorporated into the definition of $\ads\partial_{\theta}$,
corresponds to the anomalous energy shift we encountered in
(\ref{eq:energy-shift}).

Given an orthonormal basis $\{X_{i}\}$ for $\g$, the loops
$X_{i}^{n} = X_{i}z^{n}$ for $n\in\Z$ form a basis for $L\gC$
satisfying $\langle X_{i}^{n},X_{j}^{m}\rangle =
\delta_{i,j}\delta_{n,-m}$. In terms of this basis, we have
\begin{align*}
    \ads \xi
    &= -\frac{1}{4}\sum_{i,k} X_{i}^{-k} \cdot [\xi,X_{i}^{k}],
    \qquad\quad
    \ads \partial_{\theta}
    = \frac{1}{2}\sum_{j,\,k>0} i k \, X_{j}^{k}\cdot X_{j}^{-k}, \\
    \gamma
    &= -\frac{1}{24}\sum_{i,j,k,l}
    X_{i}^{-k}\cdot X_{j}^{-l}\cdot[X_{i},X_{j}]^{k+l}
    = \frac{1}{6}\sum_{i,k} X_{i}^{-k}\cdot \ads X_{i}^{k}.
\end{align*}
Note that in the expressions for $\ads\xi$ and $\gamma$, the ordering
of the factors does not matter (up to sign), since they are orthogonal
and therefore anti-commute with each other.  However, in the
expression for $\ads\partial_{\theta}$, we have
$\{X_{i}^{k},X_{i}^{-k}\} = 2$, so changing the order of the factors
shifts the operator by a constant.  Here we see normal ordering in
action, forcing us to write factors $X_{i}^{k}$ with $k$ positive on
the left and factors $X_{i}^{-k}$ with $-k$ negative on the right.  In
physics notation, this would be written as $\ads \partial_{\theta} =
-\frac{1}{4}\sum_{j,\,k\in\Z} i k \,
{\,:\,}X_{j}^{-k}\,X_{j}^{k}{\,:\,},$ where ${:\,}\xi\,\eta{\,:} =
n^{-1}( n\,\xi\wedge n\,\eta)$ denotes the normal ordered product in
the Clifford algebra.  (This colon notation is misleading as it is not
a map on the Clifford algebra but rather an instruction to replace all
Clifford products between the colons with normal ordered products.)

If the Lie algebra $\g$ is not simple, then
Theorem~\ref{th:loop-superalgebra} still holds, albeit with
slightly modified commutation relations. For a general finite
dimensional Lie algebra $\g$, the Casimir operator
$\Delta_{\ad}^{\g} = -\tfrac{1}{2}\sum_{i}(\ad X_{i})^{2}$ no
longer takes a constant value $c_{\g}$. In this case, the role of
the quadratic element $\ads\partial_{\theta}$ is played by the
2-cocycle $\omega_{\ads}$ for the projective spin representation
$\ads$, given on $\xi,\eta\in L\g$ by
\begin{equation*}
    \omega_{\ads}(\xi,\eta) := [\ads\xi,\ads\eta]-\ads\,[\xi,\eta]=
    i\,\langle\xi,\Delta_{\ad}^{\g}\partial_{\theta}\eta\rangle,
\end{equation*}
where the Casimir operator $\Delta_{\ad}^{\g}$ acts pointwise on
the loop space $L\g$. Viewing $\omega_{\ads}$ as an element of the
Clifford algebra, we have the commutator
\begin{equation*}
    [\omega_{\ads},\xi^{*}]
    = - 2\,\iota_{\xi}\,\omega_{\ads}
    = 4i\,\bigl( \Delta^{\g}_{\ad}\partial_{\theta}\xi \bigr)^{*},
\end{equation*}
so we may also view the projective cocycle as $\omega_{\ads} =
4i\ads (\Delta^{\g}\partial_{\theta})$, where $\Delta^{\g}$ is the
formal Casimir operator in the universal enveloping algebra of
$\g$. We therefore have
\begin{equation*}
    [\omega_{\ads},\ads\xi]
    = 4i\,[\ads(\Delta^{\g}\partial_{\theta}),\ads\xi]
    = 4i\ads\bigl( \Delta^{\g}_{\ad}\partial_{\theta}\xi \bigr),
\end{equation*}
and the adjoint action of $\gamma$ in
Theorem~\ref{th:loop-superalgebra} then becomes
\begin{align*}
    [\gamma,\ads\xi]
    &=\tfrac{1}{2}\,i\,\bigl(
    \Delta_{\ad}^{\g}\partial_{\theta}\xi\bigr)^{*}, \\
    \{ \gamma, \gamma \}
    &=\tfrac{1}{4}\,\omega_{\ads}-\tfrac{1}{24}\tr_{\g}\Delta_{\ad}^{\g},
\end{align*}
with the other commutation relations remaining unchanged.
Alternatively, the projective cocycle $\omega_{\ads}$ can be
viewed as the 2-form component of the Casimir operator
\begin{equation*}
    \Delta^{L\g}_{\ads}
    = -2i\ads(\Delta^{\g}\partial_{\theta}) + \Delta^{\g}_{\ads}
    = -\tfrac{1}{2}\,\omega_{\ads} +
    \tfrac{1}{8}\tr_{\g}\Delta^{\g}_{\ad}
\end{equation*}
for the spin representation $\ads$ of $L\g$, which we discuss in
Theorem~\ref{theorem:casimir} below.

\section{The Dirac operator on $L\g$}

\label{section:dirac-lg}

Following our discussion in Section~\ref{section:dirac-g}, given
any positive energy representation
$r:\Lg\rightarrow\End(\Hilbert)$, we construct a Dirac operator
\begin{equation*}
    \dirac_{r} := \hat{r} + 1 \otimes \tfrac{1}{2}\,\Omega_{L\g}
    \in \End(\Hilbert\otimes\mathcal{S}_{L\g}),
\end{equation*}
where $\hat{r}$ is the tautological $\End(\Hilbert)$-valued 1-form
on $L\g$ given by $\hat{r}(\xi) = r(\xi)$ for all $\xi\in L\g$,
and $\Omega_{L\g}$ is the fundamental 3-form given by
$\Omega_{L\g}(\xi,\eta,\zeta) =
-\frac{1}{6}\langle\xi,[\eta,\zeta]\rangle$ for $\xi,\eta,\zeta\in
L\g$. As in the previous section, we implicitly identify
$\Lambda^{*}(L\g)^{+}$ with its image in
$\End(\mathcal{S}_{L\g})$. Written in terms of a basis
$\{X_{i}^{n}\}$ of $L\g$ satisfying $\langle
X_{i}^{n},X_{j}^{m}\rangle = \delta_{i,j}\delta_{n,-m}$, this
Dirac operator is
\begin{equation*}\begin{split}
    \dirac_{r}
    &= \sum_{i,n} X_{i}^{-n}\,r(X_{i}^{n})
    - \frac{1}{12}\sum_{i,j,m,n}
        X_{i}^{-n} \cdot X_{j}^{-m} \cdot [X_{i},X_{j}]^{n+m} \\
    &= \sum_{i,n} X_{i}^{-n}\,\Bigl(
        r(X_{i}^{n}) + \frac{1}{3} \ads X_{i}^{n} \Bigr).
\end{split}\end{equation*}
Note that all of the individual factors in this expression
(anti-)commute with each other, so $\dirac_{r}$ does indeed give a
well-defined operator on the tensor product
$\Hilbert\otimes\mathcal{S}_{L\g}$, without requiring normal
ordering or dealing with any infinite constants.

In its most general form, if we take the representation $r$ to be the
canonical inclusion $r:L\g\hookrightarrow U(L\g)$ of $L\g$ into its
universal enveloping algebra $U(L\g)$, then the corresponding
universal Dirac operator is an element of the formal completion of the
non-abelian Weil algebra $\mathcal{A} = U(\Lg)\otimes\Cl(L\g)$. 
(Alternatively, we may view $\mathcal{A}$ as the universal enveloping
algebra of the super Kac-Moody algebra $\Lg_{\text{even}}\oplus
L\g_{\text{odd}}$.)  As we saw in the previous section, the product of
two such infinite formal sums does not necessarily converge.  However,
keeping in mind that we are really working with operators on Hilbert
spaces, we can indeed extend multiplication to a suitable subspace
$\mathcal{A}^{+}$ of the formal completion, which we define as the
largest subspace for which the homomorphism
$\mathcal{A}\rightarrow\End(\Hilbert\otimes\mathcal{S}_{L\g})$ extends
to $\mathcal{A}^{+}$ for any positive energy representation $\Hilbert$
of $\Lg$.  In particular, if $\mathcal{H}$ is a faithful
representation of $U(\Lg)$---we can construct such a representation by
taking the Hilbert space direct sum of countably many irreducible
positive energy representations of $\Lg$---then the homomorphism
$\mathcal{A}^{+}\hookrightarrow\End(\Hilbert\otimes\mathcal{S}_{L\g})$
induces a product structure on $\mathcal{A}^{+}$.  Fortunately, we can
perform all of our computations here using the techniques of the
previous section, working with $U(\Lg)$-valued forms on $\Lg$.

Using this extended multiplication, the square of the Dirac
operator is
\begin{equation*}
    \dirac^{2}
    = \hat{r}^{2}
    + \{ \hat{r}, \tfrac{1}{2}\,\Omega_{L\g} \}
    + \tfrac{1}{4}\,\Omega_{L\g}^{2}.
\end{equation*}
Since $\hat{r}$ is an $\End(\Hilbert)$-valued 1-form on $L\g$, its
square is a sum $\hat{r}^{2} = (\hat{r}^{2})_{0} +
(\hat{r}^{2})_{2}$ of forms of homogeneous degrees 0 and 2. For
the degree 2 component, we have $(\hat{r}^{2})_{2} =
\hat{r}\wedge\hat{r}$, and the ``curvature''
$d\hat{r}+\hat{r}\wedge\hat{r}$ of the representation $r$ is given
by
\begin{equation*}
    (d\hat{r} + \hat{r}\wedge\hat{r})(\xi,\eta)
    = \tfrac{1}{2} \bigl(\,[r(\xi),r(\eta)] - r([\xi,\eta])\,\bigr)
    = \tfrac{1}{2}\,\omega_{r}
    (\xi,\eta),
\end{equation*}
where $\omega_{r}\in\Lambda^{2}(L\g)^{+}$ is the 2-cocycle
corresponding to the projective representation $r$. If $\g$ is
simple and $I$ is the generator of the universal central extension
of $L\g$, then $\omega_{r} = 4\,r(I)\ads\partial_{\theta}$. The
degree 0 component of $\hat{r}^{2}$ is given by the following:

\begin{theorem}\label{theorem:casimir}
    The operator
    $\Delta^{L\g}_{r}:=-\tfrac{1}{2}\bigl(\hat{r}^{2}\bigr)_{0}$ is
    called the Casimir operator for the loop group $L\g$, and if $\g$
    is simple then the Casimir operator acting on the irreducible
    positive energy representation $\Hilbert_{\blambda}$ with lowest
    weight $\blambda = (m,-\lambda,h)$ is given by
    \begin{equation}\label{eq:casimir}\begin{split}
        \Delta^{L\g}_{r}
        &= - i\,( h + c_{\g})\,( r(\partial_{\theta}) -i m )
        + \Delta_{\lambda}^{\g} \\
        &=
        -i\,( h + c_{\g} )\,r(\partial_{\theta})
        + \tfrac{1}{2} \bigl( \| \blambda - \brho_{\g} \|^{2}
        - \| \brho_{\g} \|^{2} \bigr),
    \end{split}\end{equation}
    where $\brho_{\g} = ( 0, \rho_{\g}, -c_{\g})$ and
    $c_{\g} = \Delta_{\ad}^{\g}$ is
    the value of the quadratic Casimir operator of $\g$ acting on
    the adjoint representation, and the inner product is given by
    (\ref{eq:inner-product}).
\end{theorem}

\begin{proof}
In order to simply our calculations, we first note the following
identites:
\begin{equation*}\begin{split}
    [\hat{r},\ads\xi](\eta)
    &=[r(\eta),\ads\xi]-\{\hat{r},\tfrac{1}{2}[\eta,\xi]^{\ast}\}
    = r([\xi,\eta]), \\
    [\hat{r},r(\xi)](\eta)
    &= [r(\eta),r(\xi)]
    = r([\eta,\xi])
    + \langle\eta,\partial_{\theta}\xi\rangle\,r(I) \\
    &= \bigl( [\ads\xi,\hat{r}]
    + r(I) (\partial_{\theta}\xi)^{\ast}\bigr)(\eta), \text{ and}\\
    [ \ads \xi, \hat{r}^{2} ]_{0}
    &= -[\ads\xi,(\hat{r}^{2})_{2}]_{0}
    = - [ \ads \xi, d\hat{r} ]_{0}\\
    &= 2\,(\ads I)\,\hat{r}(\partial_{\theta}\xi)
    = 2\,(\ads I)\,r(\partial_{\theta}\xi).
\end{split}\end{equation*}
We now show that the commutator of $\Delta_{r}^{L\g}$ with an
element $\xi\in L\g$ is
\begin{equation*}\begin{split}
    [ \Delta_{r}^{L\g}, r(\xi) ]
    &= -\tfrac{1}{2}\,[ \hat{r}^{2}, r(\xi) ] _{0}
    = -\tfrac{1}{2}\,\bigl\{ \hat{r}, [\hat{r},r(\xi)] \bigr\}_{0}\\
    &= \tfrac{1}{2}\,\bigl\{ \hat{r}, [\hat{r},\ads\xi] \bigr\}_{0}
    -\tfrac{1}{2}\,\bigl\{
    \hat{r},r(I)(\partial_{\theta}\xi)^{\ast}\bigr\}_{0}\\
    &= - (\ads I)\,r(\partial_{\theta}\xi)
    - r(I)\,r(\partial_{\theta}\xi) \\
    &= - \bigl[ (r(I) + \ads I )\,r(\partial_{\theta}),
        r(\xi) \bigr].
\end{split}\end{equation*}
It follows that the operator $\tilde{\Delta}_{r}^{L\g} :=
\Delta_{r}^{L\g} + (r(I)+\ads I)\,r(\partial_{\theta})$ commutes
with the action of $L\g$, and therefore takes a constant value on
each irreducible representation. Acting on the minimum energy
subspace $\Hilbert_{\blambda}(m)$ of $\Hilbert_{\blambda}$, the
only terms contributing to $\Delta_{r}^{L\g}$ are those coming
from the constant loops, and thus this constant is
\begin{equation*}
    \tilde{\Delta}_{\blambda}^{L\g}
    = \Delta^{L\g}_{r} \bigr|_{\Hilbert_{\blambda}(m)}
    + i\,(h+c_{g})\,r(\partial_{\theta}) \bigr|_{\Hilbert_{\blambda}(m)}
    = \Delta^{\g}_{\lambda} - (h + c_{\g})\,m.
\end{equation*}
The desired result then follows immediately.
\end{proof}

By definition, the 0-form component of $\hat{r}^{2}$ acts as the
identity operator on $\mathcal{S}_{L\g}$. To compute the action of
$\Delta_{r}^{L\g}$, we can therefore restrict it to the minimum
energy subspace $\mathcal{S}_{L\g}(0)$ of the spin representation.
In terms of a basis $\{X_{i}^{n}\}$, we have
\begin{equation*}\begin{split}
    \Delta_{r}^{L\g}
    &= -\frac{1}{2}\,{\sum}_{i,n}r(X_{i}^{n})\, X_{i}^{-n}\,
                  {\sum}_{j,m}r(X_{j}^{m})\, X_{j}^{-m}
                  \Bigr|_{\Hilbert\otimes\mathcal{S}_{L\g}(0)} \\
    &= -\frac{1}{2}\,{\sum}_{i,j}\Bigl(
    r(X_{i})\,r(X_{j})\,X_{i}\cdot X_{j}
    + {\sum}_{n>0}r(X_{i}^{n})\,r(X_{j}^{-n})
    \,X_{i}^{-n}\cdot X_{j}^{n}\Bigr) \\
    &= \vphantom{\frac{1}{2}}\Delta_{r}^{\g} -
    {\sum}_{i,\,n>0}r(X_{i}^{n})\,r(X_{i}^{-n}),
\end{split}\end{equation*}
which is the usual definition of the Casimir operator for a loop
group. The Casimir operator can be used to define the energy
operator $r(\partial_{\theta})$ in terms of the action of $L\g$.
The constant term $\Delta^{\g}_{\lambda}$ is sometimes
incorporated into $r(\partial_{\theta})$, in which case it is
viewed as an anomalous energy shift due to the degeneracy of the
vacuum.

Returning to our computation of the square of the Dirac operator,
we note that the cross term is given by the anti-commutator $\{
\hat{r}, \tfrac{1}{2}\,\Omega_{L\g} \} = d\hat{r}$, and we obtain
\begin{equation*}\begin{split}
    \dirac^{2}
    = -2\,\Delta_{r}^{L\g} + \tfrac{1}{2}\,\omega_{\varrho}
    - \tfrac{1}{12}\tr_{\g}\Delta_{\ad}^{\g},
\end{split}\end{equation*}
where $\omega_{\varrho}$ is the 2-cocycle corresponding to the
diagonal action $\varrho = r\otimes 1 + 1\otimes\ads$ on the
tensor product $\Hilbert\otimes\mathcal{S}_{L\g}$. If $\g$ is
simple, then this 2-cocycle is $\omega_{\varrho} =
4\,\varrho(I)\ads\partial_{\theta}$. Furthermore, if
$\Hilbert_{\blambda}$ is the irreducible positive energy
representation of $L\g$ with lowest weight $\blambda =
(m,-\lambda,h)$, then using (\ref{eq:casimir}) for the Casimir
operator, we have
\begin{equation}\label{eq:loop-dirac-square}\begin{split}
    \dirac_{\blambda}^{2}
    &= 2i\,(h+c_{\g})\,\bigl(\varrho(\partial_{\theta})-im\bigr)
      - \| \lambda + \rho_{g} \|^{2} \\
    &= 2\,\varrho(I)\,\varrho(\partial_{\theta})
    - \| \blambda - \brho_{\g} \|^{2}.
\end{split}\end{equation}
Note that unlike the finite dimensional case discussed in
\S\ref{section:dirac-g}, the square of the Dirac operator for
$L\g$ does not take a constant value on each irreducible
representation. Here, the Dirac operator fails to commute with the
diagonal action $\varrho$ of $L\g$ on
$\Hilbert\otimes\mathcal{S}_{L\g}$. Since the Dirac operator
satisfies the identity $\iota_{\xi}\dirac = \varrho(\xi)$, we have
\begin{equation*}
    [\varrho(\xi),\dirac]
    = \iota_{\xi}\dirac^{2}
    = \tfrac{1}{2}\,\iota_{\xi}\,\omega_{\varrho}
    = 2\,\varrho(I)\,\iota_{\xi}\ads\partial_{\theta}
    = -\varrho(I)(\partial_{\theta}\xi)^{*},
\end{equation*}
and thus $\dirac$ commutes only with the subalgebra
$\R\oplus\g\oplus\R$ of $\R\mathop{\tilde{\oplus}}\Lg$.

If the Lie algebra $\g$ is reductive but not simple, then the
expression (\ref{eq:loop-dirac-square}) for the square of the
Dirac operator still holds provided that the central extension
satisfies $\omega_{\varrho} = 4\,\varrho(I)
\ads\partial_{\theta}$. In other words, the invariant inner
product on $\g$ must satisfy
\begin{equation*}
    [\varrho(\xi),\varrho(\eta)] - \varrho([\xi,\eta])
    = \varrho(I) \langle \xi,\partial_{\theta}\eta \rangle
\end{equation*}
for some imaginary constant $\varrho(I)$. Given any irreducible
positive energy projective representation $\mathcal{H}$ of $L\g$,
we can always choose an invariant inner product on $\g$ such that
$\mathcal{H}\otimes\mathcal{S}_{L\g}$ is a true representation of
the corresponding central extension $\Lg$. However, this choice of
inner product depends on the representation, so this approach does
not give a universal expression for the Dirac operator.

\section{The Dirac operator on $L\g/L\h$}

\label{section:dirac-lg-lh}

As in Section~\ref{section:dirac-gh}, let $\h$ be a Lie subalgebra
of $\g$, and let $\p$ denote the orthogonal complement of $\h$
with respect to the invariant inner product on $\g$. This
orthogonal decomposition extends to the loop Lie algebra $L\g =
L\h \oplus L\p$, and the Clifford algebra decomposes as
$\Cl(L\g)\cong \Cl(L\h)\otimes\Cl(L\p)$. If $\p$ is even
dimensional, as is the case when $\h$ has the same rank as $\g$,
then we can also factor the spin representation as
$\mathcal{S}_{L\g} \cong \mathcal{S}_{L\h} \otimes
\mathcal{S}_{L\p}$, where $\mathcal{S}_{L\h}$ and
$\mathcal{S}_{L\p}$ are representations of $\Lh$ of levels
$c_{\h}$ and $c_{\g}-c_{\h}$ respectively, and the action of $L\h$
on $\mathcal{S}_{L\p}$ is
\begin{align*}
    \ads_{L\p} : L\h &\rightarrow \Lambda^{2}(L\p)^{+}
    \hookrightarrow \End(\mathcal{S}_{L\p}) \\
    \zeta &\mapsto
    (\ads_{L\p}\zeta )(\xi,\eta)
    = \tfrac{1}{4}\,\langle\xi,[\zeta,\eta]\rangle
\end{align*}
for $\zeta\in L\h$ and $\xi,\eta\in L\p$.

Given any positive energy representation $r_{L\g}$ of $\Lg$ on a
Hilbert space $\Hilbert$, its restriction gives a representation
$r_{L\h}$ of $\Lh$ on $\Hilbert$. Now consider the diagonal
representation $r_{L\h}' = r_{L\h} \otimes 1 + 1 \otimes
\ads_{L\p}$ of $\Lh$ on the tensor product $\Hilbert \otimes
\mathcal{S}_{L\p}$. Using the construction of the previous
section, we build the twisted Dirac operator
\begin{equation*}
    \dirac_{L\h}' = \hat{r}_{L\h}' + \tfrac{1}{2}\,\Omega_{L\h}
    \in \End(\Hilbert\otimes\mathcal{S}_{L\p}\otimes\mathcal{S}_{L\h})
    \cong \End(\Hilbert\otimes\mathcal{S}_{L\g}).
\end{equation*}
Noting that the diagonal action $\varrho_{L\h}' = r'\otimes 1 +
1\otimes\ads_{L\h}$ on $\Hilbert\otimes\mathcal{S}_{L\g}$ is
simply the restriction of the action $\varrho_{L\g} = r\otimes 1 +
1\otimes\ads_{L\g}$ to $L\h$, we obtain the identities
\begin{align*}
    \iota_{\zeta}\dirac_{L\h}'
    &= \varrho_{L\h}'(\zeta) = \varrho_{L\g}(\zeta), \\
    [\varrho_{L\g}(\zeta),\dirac_{L\h}']
    &= \tfrac{1}{2}\,\iota_{\zeta}\,\omega_{\varrho'}^{L\h}
    = \tfrac{1}{2}\,\iota_{\zeta}\,\omega_{\varrho}^{L\g},
\end{align*}
for $\zeta\in L\h$. The difference $\dirac_{L\g/L\h} :=
\dirac_{L\g} - \dirac_{L\h}'$ is basic with respect to $L\h$, i.e.
\begin{equation*}
    \iota_{\zeta}\dirac_{L\g/L\h} = 0,
    \qquad
    [\varrho_{L\g}(\zeta),\dirac_{L\g/L\h}] = 0,
\end{equation*}
for all $\zeta\in L\h$, and thus it can be written as the
$L\h$-equivariant operator
\begin{equation}\label{eq:dirac-gh-def}
    \dirac_{L\g/L\h} = \hat{r}_{L\p} +
    \tfrac{1}{2}\,\Omega_{L\p} \in
    \End(\Hilbert\otimes\mathcal{S}_{L\p})^{L\h},
\end{equation}
where $\hat{r}_{L\p}$ is the tautological $\End(\Hilbert)$-valued
1-form on $L\p$ given by $\hat{r}(\xi) = r(\xi)$ for $\xi\in L\g$,
and $\Omega_{L\p}$ is the fundamental 3-form given by
$\Omega_{L\p}(\xi,\eta,\zeta) =
-\tfrac{1}{6}\langle\xi,[\eta,\zeta]\rangle$ for
$\xi,\eta,\zeta\in L\p$. Writing this Dirac operator in terms of a
basis $\{X_{i}^{n}\}$ of $L\p$ satisfying $\langle
X_{i}^{n},X_{j}^{m}\rangle = \delta_{i,j}\delta_{n,-m}$, we have
\begin{equation*}\begin{split}
    \dirac_{L\g/L\h}
    &= \sum_{i,n} X_{i}^{-n}\,r(X_{i}^{n})
    - \frac{1}{12}\sum_{i,j,m,n}
        X_{i}^{-n} \cdot X_{j}^{-m} \cdot
        [X_{i},X_{j}]_{\p}^{n+m},
\end{split}\end{equation*}
where $[X,Y]_{\p}$ denotes the projection of $[X,Y]$ onto $\p$.

As we saw in the finite dimensional case, the two Dirac operators
$\dirac_{L\h}'$ and $\dirac_{L\g/L\h}$ are decoupled, or in other
words they anti-commute with each other:
\begin{equation*}
    \bigl\{ \dirac_{L\h}', \dirac_{L\g/L\h} \bigr\}
    = \bigl\{ \hat{r}'_{L\h}, \dirac_{L\g/L\h} \bigr\}
    + \bigl\{ \tfrac{1}{2}\,\Omega_{L\h}, \dirac_{L\g/L\h} \bigr\}
    = 0,
\end{equation*}
where the first summand vanishes since for all $\zeta\in L\h$ we
have
\begin{equation*}
    \bigl\{ \hat{r}'_{L\h}, \dirac_{L\g/L\h} \bigr\} (\zeta)
    = [ r'(\zeta), \dirac_{L\g/L\h} ]
    = 0,
\end{equation*}
and the second summand vanishes as the odd operators
$\tfrac{1}{2}\,\Omega_{L\h}$ and $\dirac_{L\g/L\h}$ act on
distinct representations $\mathcal{S}_{L\h}$ and
$\Hilbert\otimes\mathcal{S}_{L\p}$ and therefore anti-commute.
Since these two operators are decoupled, the square of the Dirac
operator on $L\g/L\h$ is
\begin{equation*}\begin{split}
    \dirac_{L\g/L\h}^{2}
    &= \bigl( \dirac_{L\g} \bigr)^{2}
    - \bigl( \dirac_{L\h}' \bigr)^{2} \\
    &= -2\,\bigl( \Delta_{r}^{L\g} - \Delta_{r'}^{L\h} \bigr)
     + \tfrac{1}{2}\,\bigl(\omega_{\varrho}^{L\g} -
     \omega_{\varrho'}^{L\h} \bigr)
     + \tfrac{1}{12}\,
     \bigl( \tr_{\g}\Delta_{\ad}^{\g} - \tr_{\h}\Delta_{\ad}^{\h}
     \bigr).
\end{split}\end{equation*}
Now consider the case where $\g$ is simple, $\h$ is reductive, and
$\Hilbert_{\blambda}$ is the irreducible positive energy
representation of $L\g$ with lowest weight $\blambda$. Since
$\dirac_{L\g/L\h}$ is an $L\h$-equivariant operator on
$\Hilbert_{\blambda}\otimes\mathcal{S}_{L\p}$, it is a constant on
each of the irreducible subrepresentations of $L\h$. If
$\mathcal{U}_{\bmu}$ is the irreducible positive energy
representation of $\Lh$ with lowest weight $\bmu$, then using
(\ref{eq:loop-dirac-square}), we see that the square of the Dirac
operator takes the value
\begin{equation}\label{eq:lg-lh-dirac-square}\begin{split}
    (\dirac_{L\g/L\h})^{2}\bigr|_{\bmu}
    &= 2\,\varrho_{L\g}(I)\,\varrho_{L\g}(\partial_{\theta})
    - \| \blambda - \brho_{\g} \|^{2}\\
    &- 2\,\varrho_{L\h}'(I)\,\varrho_{L\h}'(\partial_{\theta})
    + \| \bmu - \brho_{\h} \|^{2}
    = - \|\blambda-\brho_{\g}\|^{2}+\|\bmu-\brho_{\h}\|^{2},
\end{split}\end{equation}
on the $\mathcal{U}_{\bmu}$ components of
$\Hilbert_{\blambda}\otimes\mathcal{S}_{L\p}$. Note that the
non-constant terms vanish since $\varrho_{L\g}$ and
$\varrho_{L\h}'$ agree on $\R\mathop{\tilde{\oplus}}\Lh$.

Note that in the above construction, we are using the invariant
inner product on $\h$ obtained by restricting our invariant inner
product on $\g$. When $\g$ is simple, we use the basic inner
product on $\g$, which is normalized so that
$\|\alpha_{\text{max}}\|^{2} = 2$, where $\alpha_{\text{max}}$ is
the highest root of $\g$. We recall that the basic inner product
corresponds to the universal central extension $\Lg$ of $L\g$,
which in turn restricts to give a (not necessarily universal)
central extension $\Lh$ of $L\h$. Nevertheless, given any positive
energy representation $\Hilbert$ of $\Lg$, the tensor product
$\Hilbert\otimes\mathcal{S}_{L\g}$ is a true representation of
this central extension $\Lh$. So, if $\h$ is reductive, then the
squares of the Dirac operators $\dirac_{L\h}'$ and
$\dirac_{L\g/L\h}$ are indeed of the form given by
(\ref{eq:loop-dirac-square}) and (\ref{eq:lg-lh-dirac-square}).

On the other hand, if $\g$ is not simple but rather semi-simple,
then the basic inner product on $\g$ is normalized so that
$\|\alpha_{i}\|^{2} = 2$, where the $\alpha_{i}$ are the highest
roots of each of the simple components of $\g$. In this case, a
projective positive energy representation of $L\g$ is not
necessarily a true representation of the corresponding central
extension $\Lg$, so the expression (\ref{eq:loop-dirac-square})
for the square of the Dirac operator on $L\g$ is not universal.
However, if $\h$ is reductive, then the expression
(\ref{eq:lg-lh-dirac-square}) for the square of the Dirac operator
on $L\g/L\h$ does still hold, as the non-constant terms must
vanish since the operator commutes with the action of $L\h$.

\section{The kernel of the Dirac operator}

\label{section:kernel}

Given a linear operator $d:V\rightarrow W$ between two finite
dimensional vector spaces, the alternating sum of the dimensions
in the exact sequence
\begin{equation*}
    0 \longrightarrow \Ker d \longrightarrow V
    \stackrel{d}{\longrightarrow}
    W \longrightarrow \Coker d \longrightarrow 0
\end{equation*}
vanishes, and it follows that $\Index d = \dim V - \dim W$.
Furthermore, if $V$ and $W$ are $G$-modules and the operator $d$
is $G$-equivariant, then the analogous result $\Index_{G}d = V -
W$ holds in the representation ring $R(G)$. In the infinite
dimensional case, this result does not necessarily hold, but for
representations of loop groups, it does hold provided that the
representations are of finite type and that the operator commutes
with rotating the loops.

\begin{lemma}\label{lemma:index}
    If $\mathcal{V}$ and $\mathcal{W}$ are representations of $LG$ of
    finite type, and $\mathcal{D}:\mathcal{V}\rightarrow\mathcal{W}$
    is an $S^{1}\semi LG$-equivariant linear operator, then its
    $LG$-equivariant index is the virtual representation
    $\Index_{LG}\mathcal{D}=\mathcal{V}-\mathcal{W}$.
\end{lemma}

\begin{proof}
    Since $\mathcal{D}$ is $S^{1}$-equivariant, it respects the
    decompositions of $\mathcal{V}$ and $\mathcal{W}$ into their
    constant energy subspaces, and it can be written in the block
    diagonal form $\mathcal{D} = \bigoplus_{k\in\Z}\mathcal{D}_{k}$,
    with $\mathcal{D}_{k}:\mathcal{V}(k)\rightarrow \mathcal{W}(k)$. 
    If both $\mathcal{V}$ and $\mathcal{W}$ are of finite type, then
    each of the subspaces $\mathcal{V}(k)$ and $\mathcal{W}(k)$ is a
    finite dimensional $G$-module, and so the $S^{1}\times
    G$-equivariant index of $\mathcal{D}$ is given by the
    $R(G)$-valued formal power series
    \begin{equation*}
    \Index_{S^{1}\times G}\mathcal{D}
    =\sum_{k\in\Z}z^{k}\,\bigl( \mathcal{V}(k) - \mathcal{W}(k) \bigr)
    =\sum_{k\in\Z}z^{k}\,\mathcal{V}(k)-\sum_{k\in\Z}z^{k}\,\mathcal{W}(k).
    \end{equation*}
    Since a representation of the full loop group $LG$ is uniquely
    determined by its constant energy components, the
    $LG$-equivariant index must therefore be the difference of the
    domain and the range, hence
    $\Index_{LG}\mathcal{D} = \mathcal{V} - \mathcal{W}$.
\end{proof}

Returning to the notation of the previous section, let $\g$ be
semi-simple, and let $\h$ be a reductive subalgebra of $\g$ with
maximal rank. If we decompose the spin representation as
$\mathcal{S}_{L\p} =
\mathcal{S}_{L\p}^{+}\oplus\mathcal{S}_{L\p}^{-}$, the Dirac
operator $\dirac_{L\g/L\h}$ interchanges the two half-spin
representations and can thus be written as the sum of the
operators
\begin{align*}
    \dirac_{L\g/L\h}^{+} &:
    \Hilbert_{\blambda}\otimes\mathcal{S}_{L\p}^{+}
    \rightarrow
    \Hilbert_{\blambda}\otimes\mathcal{S}_{L\p}^{-}, \\
    \dirac_{L\g/L\h}^{-} &:
    \Hilbert_{\blambda}\otimes\mathcal{S}_{L\p}^{-}
    \rightarrow
    \Hilbert_{\blambda}\otimes\mathcal{S}_{L\p}^{+},
\end{align*}
where $\dirac_{L\g/L\h}^{-}$ is the adjoint of
$\dirac_{L\g/L\h}^{+}$. When we introduced $\dirac_{L\g/L\h}$ in
(\ref{eq:dirac-gh-def}), we showed that it is $L\h$-equivariant,
and all of our Dirac operators clearly commute with the generator
$\partial_{\theta}$ of rotations of the loops. The operator
$\dirac_{L\g/L\h}$ is therefore $S^{1}\semi LH$-equivariant, and
since its domain and range are both of finite type, we may apply
Lemma~\ref{lemma:index}. The $LH$-equivariant index of
$\dirac_{L\g/L\h}^{+}$ is thus the difference
\begin{equation*}
    \Ker\dirac_{L\g/L\h}^{+} - \Ker\dirac_{L\g/L\h}^{-}
    = \Hilbert_{\blambda}\otimes\mathcal{S}_{L\p}^{+}
    - \Hilbert_{\blambda}\otimes\mathcal{S}_{L\p}^{-},
\end{equation*}
which is given by the homogeneous Weyl-Kac formula
(\ref{eq:homogeneous-weyl-kac}).

On the other hand, to compute the kernel of $\dirac_{L\g/L\h} =
\dirac_{L\g/L\h}^{+}\oplus\dirac_{L\g/L\h}^{-}$ we proceed as in
the computation of the kernel of the finite dimensional operator
$\dirac_{\g/\h}$ in \S\ref{section:dirac-gh}. In fact, the proofs
of Lemmas \ref{lemma:1} and \ref{lemma:2} apply equally well in
the Kac-Moody setting using the decomposition
    $\mathcal{S}_{L\g/\mathfrak{t}} \cong \mathcal{S}_{L\p} \otimes
    \mathcal{S}_{L\h/\mathfrak{t}},$
and we obtain

\begin{lemma}
For each $c\in \mathcal{C}$, the irreducible representation
$\mathcal{U}_{c\bullet\blambda}$ of $\Lh$ with lowest weight
$c\bullet\blambda = c(\blambda-\brho_{\g}) + \brho_{\h}$ occurs
exactly once in the decomposition of
$\Hilbert_{\blambda}\otimes\mathcal{S}_{L\p}$.
\end{lemma}

\begin{lemma}
If $\bmu$ is a weight of
$\Hilbert_{\blambda}\otimes\mathcal{S}_{L\p}$ satisfying $\|
\bmu-\brho_{\h} \|^{2} = \| \blambda-\brho_{\g} \|^{2}$, then
there exists a unique affine Weyl element $w\in\Waff_{\g}$ such
that $\bmu-\brho_{\h} = w(\blambda-\brho_{\g})$.
\end{lemma}

Then, in light of our formula (\ref{eq:lg-lh-dirac-square}) for
the square of the Dirac operator $\dirac_{L\g/L\h}$, we
immediately obtain the loop group analogue of
Theorem~\ref{theorem:dirac-gh-kernel}.

\begin{theorem}\label{theorem:dirac-lgh-kernel}
    Let $\g$ be a semi-simple Lie algebra with a maximal rank
    reductive Lie subalgebra $\h$.  Let $\Hilbert_{\blambda}$ and
    $\mathcal{U}_{\bmu}$ be the irreducible representations of $\Lg$
    and $\Lh$ with lowest weights $\blambda$ and $\bmu$.  The kernel
    of the operator $\dirac_{L\g/L\h}$ on
    $\Hilbert_{\blambda}\otimes\mathcal{S}_{L\p}$ is
    \begin{equation*}
        \Ker\dirac_{L\g/L\h}
        ={\bigoplus}_{c\in\mathcal{C}}\mathcal{U}_{c\bullet\blambda},
    \end{equation*}
    where $c\bullet\blambda = c(\blambda-\brho_{\g})+\brho_{\h}$, and
    $\mathcal{C}\subset\Waff_{\g}$ is the subset of affine Weyl
    elements which map the fundamental Weyl alcove for $\g$ into the
    fundamental alcove for $\h$.
\end{theorem}

Comparing this result to the homogeneous Weyl-Kac formula
(\ref{eq:homogeneous-weyl-kac}), we obtain
\begin{equation*}
    \Ker\dirac_{L\g/L\h}^{+} = \bigoplus_{(-1)^{c} = +1}
    \mathcal{U}_{c\bullet\blambda}, \qquad
    \Ker\dirac_{L\g/L\h}^{-} = \bigoplus_{(-1)^{c} = -1}
    \mathcal{U}_{c\bullet\blambda}.
\end{equation*}
Taking the kernels of these Dirac operators therefore gives an
explicit construction for the multiplet of signed representations
of $\Lh$ corresponding to any given irreducible positive energy
representation of $\Lg$.

\end{document}